\documentclass[final,onefignum,onetabnum]{siamart220329}

\usepackage{lipsum}
\usepackage{amsfonts}
\usepackage{graphicx}
\usepackage{epstopdf}
\usepackage{algorithmic}
\usepackage{amssymb}
\usepackage{booktabs}
\usepackage{svg}
\usepackage{setspace}
\ifpdf
  \DeclareGraphicsExtensions{.eps,.pdf,.png,.jpg}
\else
  \DeclareGraphicsExtensions{.eps}
\fi

\graphicspath{{figs/}}

\newsiamremark{remark}{Remark}
\newsiamremark{hypothesis}{Hypothesis}
\crefname{hypothesis}{Hypothesis}{Hypotheses}
\newsiamthm{claim}{Claim}

\headers{Opt. Cheby. Sm. and V-cycle}{M. Phillips and P. Fischer}

\title{Optimal Chebyshev Smoothers and One-sided V-cycles\thanks{This work was funded by the Exascale Computing Project under contract no. 17-SC-20-SC}}

\author{Malachi Phillips\thanks{Department of Computer Science, University of Illinois at Urbana-Champaign, Urbana IL 61801
  (\email{malachi2@illinois.edu}).
  }
\and Paul Fischer$^{\dagger}$\thanks{Department of Mechanical Science and Engineering,
University of Illinois at Urbana-Champaign, Urbana IL 61801
(\email{fischerp@illinois.edu}).
  }
}

\usepackage{amsopn}

\newcommand{\norm}[1]{{\left\| {#1}\right\|}}
\newcommand{\isdef}{\mathrel{\mathop:}=}

\newcommand{\iter}{k}
\newcommand{\nnz}[1]{\operatorname{nnz}({#1})}

\DeclareMathOperator*{\argmax}{arg\,max}
\DeclareMathOperator*{\argmin}{arg\,min}

\def\dO{\partial \Omega}

\def\dt{ \Delta t }

\def\cH{{\cal H}}

\def\scriptO{{{\it O}\kern -.42em {\it `}\kern + .20em}}
\def\RR{{{\rm l}\kern - .15em {\rm R} }}
\def\PP{{{\rm l}\kern - .15em {\rm P} }}
\def\L2{{{\sf L}^2}}
\def\H1{{{\sf H}^1}}
\def\PN2{{\PP_{N}-\PP_{N-2}}}

\def\complex{{{\rm C} \kern - .53em {\rm l} \kern + .38em}}

\def\a1{{ | \lambda_{\min} |}}

\def\l1{{   \lambda_{\min}  }}

\def\tlam {{\tilde \lambda}}

\def\bu0{{\underline {\bf 0}}}

\def\bu{{\bf u}}

\def\bx{{\bf x}}

\def\Oh{{\hat \Omega}}

\def\ub{{\underline b}}

\def\uf{{\underline f}}

\def\ur{{\underline r}}
\def\us{{\underline s}}

\def\uu{{\underline u}}

\def\u0{{\underline 0}}
\def\1u{{\underline 1}}

\def\uub{{\bar {\underline u}}}

\ifpdf
\hypersetup{
  pdftitle={Optimal Chebyshev Smoothers and One-sided V-cycles},
  pdfauthor={M. Phillips and P. Fischer}
}
\fi

\begin{document}

\maketitle

\begin{abstract}
The solution to the Poisson equation arising from the spectral element discretization
of the incompressible Navier-Stokes equation requires robust preconditioning strategies.
One such strategy is multigrid.
To realize the potential of multigrid methods, effective smoothing strategies are needed.
Chebyshev polynomial smoothing proves to be an effective smoother.
However, there are several improvements to be made, especially at the cost of symmetry.
For the same cost per iteration, a symmetric V-cycle with $\iter$ order Chebyshev polynomial smoothing
may be substituted with a one-sided V-cycle with order $2\iter$ Chebyshev polynomial smoothing,
wherein the smoother is omitted on the up-leg of the V-cycle.
The choice of omitting the post-smoother in favor of higher order Chebyshev pre-smoothing
is shown to be advantageous in cases where the multigrid approximation property constant, $C$, is large.
Results utilizing Lottes's fourth-kind Chebyshev polynomial smoother are shown.
These methods demonstrate substantial improvement over the standard Chebyshev polynomial smoother.
The authors demonstrate the effectiveness of this scheme in $p$-geometric multigrid,
as well as a 2D model problem with finite differences.
\end{abstract}



\section{Introduction}
\label{sec:introduction}

Chebyshev smoothing was introduced in the context of parallel multigrid
(MG) methods in \cite{adams2003parallel},  where it was established that Chebyshev
smoothing was competitive with Gauss-Seidel smoothing even in serial 
computing applications. Here, we explore several variations on Chebyshev
smoothing for the Poisson problem in general domains. Our primary aim is to
develop fast highly-scalable solvers for the pressure sub-step in time
advancement of the Navier-Stokes (NS) equations, particularly for discretizations
based on the spectral element method (SEM). Many of the finding, however,
would apply in more general settings.

Our target problem is to solve a sequence of Poisson problems,
\begin{align}
  -\nabla^2 \tilde u = \tilde f \text{ for } \tilde u,\tilde f \in \Omega \subset \RR^d \mapsto \RR.
  \label{eqn:poisson}
\end{align}
The weak formulation is written as:
{\em find $u^m(\bx) \in X^N_0 \subset \cH^1_0$ such that}
\begin{eqnarray}
\int_{\Omega} \nabla v \cdot \nabla u \, dV &=& \int_{\Omega} \, v \, f^m \, dV
\hspace*{.1in} \forall \, v \, \in \, X^N_0,
\end{eqnarray}
where $f^m(\bx)$ is the data and $u^m(\bx)$ is the corresponding solution
field at some time instant $t^m$, $m=1,$ 2, $\dots$
Here, $\Omega \subset \RR^d$ is the computational domain in $d$ (=1, 2, or 3)
space dimensions; $\cH^1_0(\Omega)$ is the standard Sobolev space
comprising functions that vanish on a subset of the boundary,
$\dO_D \subset \dO$, 
are square-integrable on $\Omega$, and  whose
gradient is also square-integrable; and $X_0^N = \operatorname{span}\{ \phi_j(\bx) \}$
is the finite-dimensional trial/test space associated with a Galerkin
formulation of the Poisson problem. The discrete problem statement
is expressed as $A \uu^m = \ub^m$, where $\uu^m$ is the vector of
basis coefficients at $t^m$ and $A$ is the symmetric-positive
definite (SPD) matrix with
\begin{eqnarray}
a_{ij} &=& \int_{\Omega} \, \nabla \phi_i \cdot \nabla \phi_j \, dV.
\end{eqnarray}

We point out that the need to solve a sequence of problems differs from solving
a single problem in several significant ways. First, solver set-up costs are
typically amortized over thousands of right-hand sides and are therefore
largely irrelevant to our cost concerns.  Second, the solution is typically
devoid of significant low wave-number content because we solve only for a
perturbed solution, $\delta \uu^m := \uu^m - \uub$, where $\uub$ is an
initial guess.  If we take $\uub=\uu^{m-1}$ then the initial residual $\ur_0
= \ub - A \uu^{m-1} = O(\dt)$. This result is improved to $O(\dt^l)$ by
projecting $\uu^m$ onto the space of prior solutions, $\{ \uu^{m-1} \; \dots \;
\uu^{m-l}  \}$ \cite{chanwan97,fischer1998projection}. Finally, with an initially small
residual, GMRES is likely to converge in just a few iterations, which obviates
the need for restarts and mitigates the $O(k^2)$ complexity terms in a
$k$-iteration GMRES solve. This latter observation puts less pressure on
requiring a symmetric preconditioner since one can retain the full benefits of
using Krylov subspace projection (KSP) without resorting to conjugate gradient
iteration. With these circumstances in mind, we will drop the superscript $m$
in the sequel.

We note that Chebyshev smoothers have gained a lot of attention recently.
Kronbichler and co-workers \cite{kronbichler2019multigrid,fehn18b, fehn20} have employed Chebyshev smoothing
for discontinuous Galerkin discretizations of the NS equations.
Rudi and coworkers employ algebraic multigrid (AMG) with Chebyshev smoothing \cite{rudi2015extreme}.
Similarly, Chebyshev smoothing is considered by Sundar and coworkers
as a multigrid smoother for high-order continuous finite element discretizations \cite{sundar2015comparison}.

A major difference here is that we consider Chebyshev in conjunction with
additive Schwarz methods (ASM) \cite{wid88, fischer97,lottes2005hybrid,phillips2022tuning}
and restrictive additive Schwarz (RAS) \cite{cai1999restricted} in place of point-Jacobi
smoothing.
The principal idea is to use ASM or RAS to eliminate high wave number content.
In the case of the spectral element method, local Schwarz
solves can be effected at a cost that is comparable to forward operator
evaluation through the use of fast diagonalization
\cite{lynch64,fischer00a,lottes2005hybrid}.
Another critical aspect of the current context is that many of our applications
are targeting exascale platforms and beyond, where compute is performed  on
tens of thousands of GPUs for which the relative cost of global communication and
hence, coarse-grid solves, is high \cite{min2022optimization}. In such cases, it
often pays to have high-quality and broad bandwidth smoothing, such as provided
by Chebyshev, in order to reduce the number of visits to the bottom of the V-cycle
where the expensive coarse-grid solve is invoked.

Here, we explore a seemingly simple question:
{\em Given $2k$ smoothing iterations,
what is the optimal choice of $m$ pre-smoothing
and $n$ post-smoothing applications, where $m+n=2k$?
}
More specifically, in the Chebyshev context, the question is
what order $m$ pre-smoothing and order $n$ post-smoothing should be used
at the same cost per iteration.
An additional and important point to this question is, {\em What kind of Chebyshev
smoothing should be used?}   One could use standard 1st-kind Chebyshev
polynomials with tuned parameters. (Recall, we can afford significant tuning
overhead.)  Or, one could use standard or optimized 4th-kind Chebyshev
polynomials that were proposed in recent work by Lottes. (See \cite{lottes2022optimal}
and references therein.)   We explore these questions under several different
conditions: using finite differences and spectral elements discretizations
and using Jacobi, ASM, or RAS as the basic smoother.

The structure of this paper is as follows.
\Cref{sec:chebyshev-mg} outlines the multigrid V-cycle and
Chebyshev smoothers.
2D finite difference Poisson results on varying aspect ratio grids, along with a comparison
between theoretical and observed multigrid error contraction rates, are
presented in \cref{sec:fd-gmg}.
Spectral element (SE)-based pressure Poisson preconditioning schemes
implemented in the
scalable open-source CFD code, nekRS \cite{fischer2022nekrs},
are presented in \cref{sec:pmg-preconditioners}.
nekRS started as a fork of libParanumal \cite{libp} and uses
highly optimized kernels based on the Open Concurrent Compute Abstraction (OCCA) \cite{medina2014occa}.
Special focus is given to performance on large-scale GPU-based
platforms such as OLCF's Summit.
Cases for the stationary Poisson and pressure Poisson problem
arising from the NS equations are shown in \cref{sec:cases}.
Results for the cases in \cref{sec:cases} are shown in \cref{sec:results}.
\Cref{sec:conclusions} concludes the paper.

\section{Multigrid and Chebyshev Smoothers}
\label{sec:chebyshev-mg}
Let us consider the V-cycle algorithm to solve the SPD matrix $A$.
Suppose that levels $j=0,\dots,\ell$ are used in the V-cycle, with $A=A_0$.
Let us denote the interpolation operator mapping entries from grid $j+1$ to $j$ by $P_{j+1}^j$ for $j=0,\dots,\ell-1$.
The sequence of matrices corresponding to each level are typically constructed in a Galerkin fashion, with
\begin{equation}\label{eqn:galerkin-coarse-grid}
  A_{j+1} = \left(P_{j+1}^j\right)^T A_j P_{j+1}^j, j=0,\dots,\ell-1.
\end{equation}
The multiplicative error propagator for a single V-cycle is given recursively by
\begin{align}\label{eqn:recursive-v-cycle-err-prop}
  E_j &= I - M_j^{-1} A_j \nonumber \\
      &= {G'_j}\left(I - P_{j+1}^j M_{j+1}^{-1}\left(P_{j+1}^j\right)^T A_j\right){G_j}, j=0,\dots,\ell-1,
\end{align}
with $M_{\ell}^{-1} \isdef A_{\ell}^{-1}$.
${G_j}$ and ${G'_j}$ are the smoother iteration matrices for the pre- and post-smoothing
iteration matrices, respectively.
For example, $G_j=(I - \omega S_j A_j)^k$ corresponds to $k$ steps of the simple smoothing iteration
\begin{equation}\label{eqn:simple-smoothing-iteration}
  \left(\underline x_{i+1}\right)_j= \left(\underline x_{i}\right)_j+\omega S_j(\underline b_j - A_j \left(\underline x_i\right)_j),
\end{equation}
where $S_j$ is the smoother for the $j$th level, such as Jacobi with $S_j = \operatorname{diag}(A_j)^{-1}A_j$.
$G'_j$ and $G_j$ are not necessarily the same -- in fact, we will consider
the choice of $G'_j = G_j$ (symmetric post-smoothing) as well as $G'_j = I$ (omitting post-smoothing),
among others.

For later use, we will need to consider the one-sided V-cycle, which is sufficient
for the analysis of general V-cycles \cite{mccormick1985multigrid,lottes2022optimal}.
Similar as before, the ``fine-to-coarse''
\begin{align}\label{eqn:recursive-ftc-err}
  \left(E_{\searrow}\right)_j &= I - \left(M_{\searrow}\right)_j^{-1} A_j \nonumber \\
      &= \left(I - P_{j+1}^j \left(M_{\searrow}\right)_{j+1}^{-1}\left(P_{j+1}^j\right)^T A_j\right){G_j}, j=0,\dots,\ell-1,
\end{align}
and ``coarse-to-fine''
\begin{align}\label{eqn:recursive-ctf-err}
  \left(E_{\nearrow}\right)_j &= I - \left(M_{\nearrow}\right)_j^{-1} A_j \nonumber \\
      &= {G'_j}\left(I - P_{j+1}^j \left(M_{\nearrow}\right)_{j+1}^{-1}\left(P_{j+1}^j\right)^T A_j\right), j=0,\dots,\ell-1,
\end{align}
error propagators are defined, with
$\left(M_{\searrow}\right)_{\ell}^{-1} = \left(M_{\nearrow}\right)_{\ell}^{-1} = A_{\ell}^{-1}$.
Let $E_{\searrow}=\left(E_{\searrow}\right)_0$ and $E_{\nearrow} = \left(E_{\nearrow}\right)_0$.
The general V-cycle, therefore, is the product of the ``fine-to-coarse'' and ``coarse-to-fine'' error propagators,
\begin{equation}\label{eqn:general-v-cycle}
  E_V = E_{\nearrow}E_{\searrow}.
\end{equation}

Given $k_j$ iterations of the smoothing iteration in \cref{eqn:simple-smoothing-iteration} for level $j$,
is it possible to construct a better order $k_j$ polynomial than $p_{k_j}(S_jA_j) = G_{k_j} =(I - \omega S_j A_j)^{k_j}$?
Following \cite{saad2003iterative,adams2003parallel}, we wish to solve:
\begin{equation}\label{eqn:opt-polynomial-cheby}
  \min_{p_k\in\mathbb P_k, p_k(0)=1} \max_{\lambda \in [\lambda_{min}, \lambda_{max}]} |p(t)|.
\end{equation}
where $\mathbb P_k$ is the space of all polynomials with degree less than or equal to $k$.
Taking $\cal E$ to be the interval $[\lambda_{min}, \lambda_{max}]$,
the solution to the \emph{minimax} problem in \cref{eqn:opt-polynomial-cheby} are the shifted and scaled Chebyshev polynomials
of the 1st-kind
\begin{equation}\label{eqn:1st-kind-cheby-polynomials}
  \hat{T_k}(\lambda) = \dfrac{1}{\sigma_k} T_k\left(\dfrac{\theta-\lambda}{\delta}\right) \text{ with } \sigma_k\isdef T_k\left(\dfrac{\theta}{\delta}\right).
\end{equation}
$T_k(\cdot)$ is the Chebyshev polynomial of the 1st-kind of order $k$;
$\theta$ is the midpoint of the interval $[\lambda_{min}, \lambda_{max}]$,
\begin{equation*}
  \theta = \dfrac{\lambda_{min}+\lambda_{max}}{2};
\end{equation*}
and $\delta$ is the mid-width of the interval,
\begin{equation*}
  \delta = \dfrac{\lambda_{max}-\lambda_{min}}{2}.
\end{equation*}
The Chebyshev polynomials of the 1st-kind enjoy a three-term recurrence relation that is used
to derive \cref{alg:cheby} \cite{saad2003iterative}.
While $\lambda_{max}$ is taken to be the largest eigenvalue of $S_jA_j$,
how should $\lambda_{min}$ be chosen?
$\lambda_{min}$ is chosen as a small factor of $\lambda_{max}$.
Previous works considered factors such as 1/30 \cite{adams2003parallel},
3/10 \cite{baker2011multigrid}, 1/4 \cite{sundar2015comparison}, and 1/6 \cite{zhukov2015multigrid}.
Two approaches utilizing the 1st-kind Chebyshev smoother are considered in this study.
The first uses $\lambda_{min} = 0.1\lambda_{max}$ and is denoted as $1^{st}$-Cheb.
The second optimizes the for $\lambda_{min}$ and is denoted as $1^{st}$-Cheb, $\lambda_{min}^{opt}$.

\begin{algorithm}\small
  \caption{\small Chebyshev smoother, 1st-kind}
  \label{alg:cheby}
  \begin{algorithmic}
    \STATE{$\theta = \dfrac 1 2 (\lambda_{max}+\lambda_{min})$, $\delta= \dfrac 1 2 (\lambda_{max}-\lambda_{min})$, $\sigma = \dfrac \theta \delta$, $\rho_0 = \dfrac 1 \sigma$}
    \STATE{$  \underline x_0 = \underline x, \underline r_0 =   S(  \underline b -   A   \underline x_0)$, $  \underline d_0 = \dfrac 1 \theta  \underline r_0$} \label{alg:cheby:init-res-eval}
    \FOR{$i=1,\dots,\iter-1$}
      \STATE{$  \underline x_{i} =   \underline x_{i-1} +   \underline d_{i-1}$}
      \STATE{$  \underline r_{i} =   \underline r_{i-1} -   S   A   \underline d_{i-1}$, $\rho_{i} = \dfrac{1}{2\sigma-\rho_{i-1}}$}
      \STATE{$  \underline d_{i}=\rho_{i}\rho_{i-1}   \underline d_{i-1} + \dfrac {2\rho_{i}}{\delta}   \underline r_{i}$}
    \ENDFOR
    \STATE{$  \underline x_{\iter} =   \underline x_{\iter-1} +   \underline d_{\iter-1}$}
    \RETURN{$  \underline x_{\iter}$}
  \end{algorithmic}
\end{algorithm}

\begin{figure}
  \centering
  \def\svgwidth{\textwidth}
  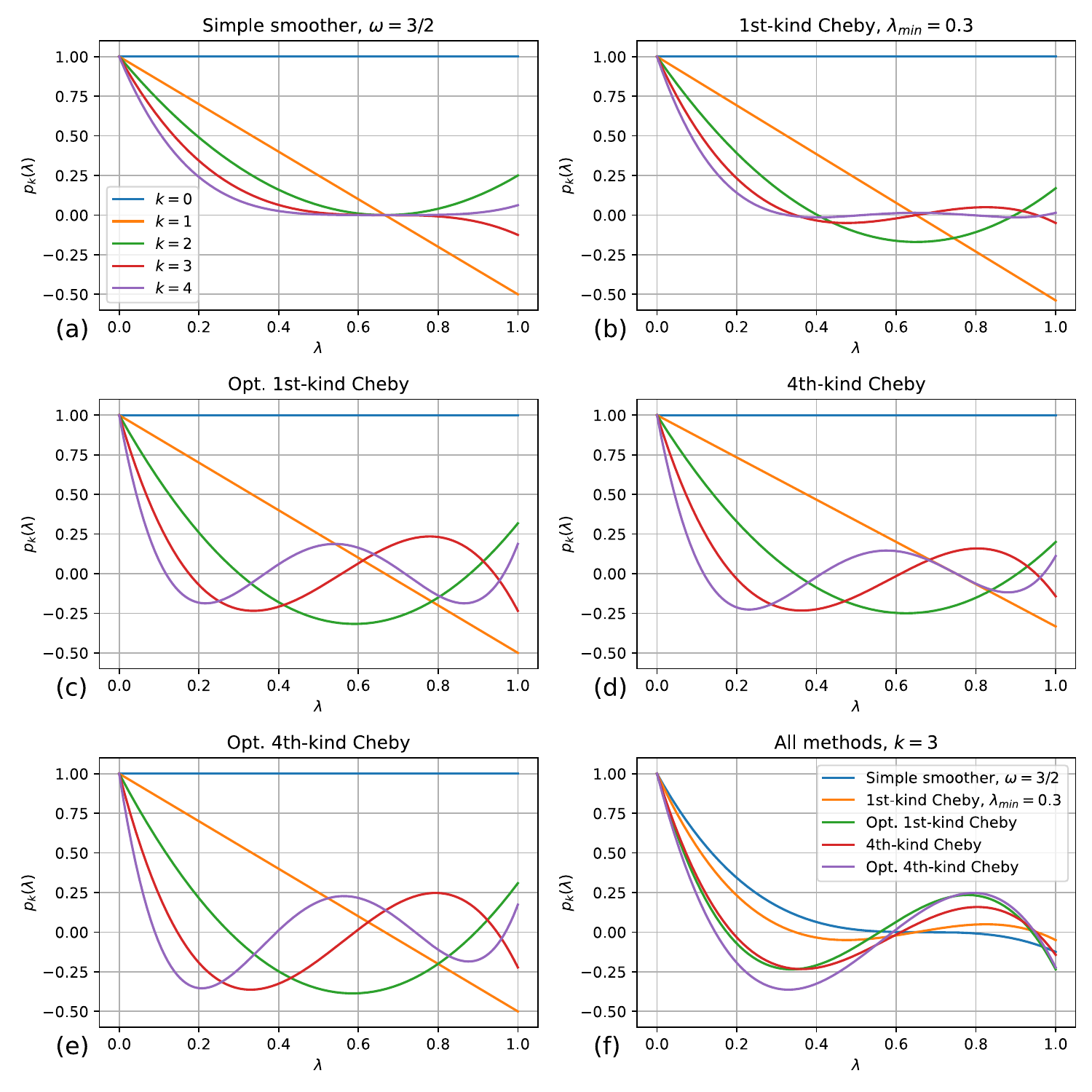
  \caption{
    \label{fig:poly-smoothers}
    $p_k(\lambda)$ for various polynomial smoothers with varying polynomial orders, $k$.
    \cref{fig:poly-smoothers}a simple polynomial smoothing with $p_k(\lambda) = (1 - 3/2 \lambda)^k$.
    \cref{fig:poly-smoothers}b smoothing with Chebyshev polynomials of the 1st-kind with $\lambda_{min}=0.3$.
    \cref{fig:poly-smoothers}c smoothing with Chebyshev polynomials of the 1st-kind, $\lambda_{min}$ chosen to minimize \cref{eq:gamma-defn}.
    \cref{fig:poly-smoothers}d smoothing with Chebyshev polynomials of the 4th-kind, $\beta_i = 1$.
    \cref{fig:poly-smoothers}e smoothing with Chebyshev polynomials of the 4th-kind, $\beta_i$ chosen to minimize \cref{eq:gamma-defn}.
    \cref{fig:poly-smoothers}f all polynomials considered, with $k=3$.
  }
\end{figure}

Is it possible to further improve the polynomial smoother and remove the ad-hoc $\lambda_{min}$ parameter?
The polynomial smoother can be chosen in such a way
to minimize an error bound \cite{lottes2022optimal}, such as the \emph{two-level} bound proposed by Hackbusch \cite{hackbusch1982multi}
in \cref{eqn:hackbusch-bound}.
Without loss of generality, let $\rho(SA)=1$.
Let $G=G_k(SA)$ be a $k$-order polynomial in $SA$.
The \emph{two-level} bound from Hackbusch \cite{hackbusch1982multi}
is re-written as \cite{lottes2022optimal}:
\begin{align}\label{eqn:hackbusch-bound}
  \norm{E_{\searrow}}_A = \norm{(I-PA_c^{-1}P^TA)G_k}_A \nonumber \\
                    \leq C_0^{1/2} \sup_{0 < \lambda \leq 1} \lambda^{1/2}|p_k(\lambda)|,
\end{align}
where $C_0$ is the multigrid approximation property constant, which for a given level $j$ is
\begin{align}\label{eq:C-def}
  C_j &\isdef \norm{A_j^{-1} - P_{j+1}^jA_{j+1}^{-1}\left(P_{j+1}^j\right)^T}_{A_j,S_j}^2 \nonumber \\
      &\isdef \sup_{\norm{\underline f}_{S_j}\leq 1} \norm{\left(A_j^{-1}-P_{j+1}^jA_{j+1}^{-1}\left(P_{j+1}^j\right)^T\right)\underline f}_{A_j}^2.
\end{align}
Solving the \emph{weighted} minimax problem in \cref{eqn:hackbusch-bound} yields the solution \cite{lottes2022optimal}:
\begin{equation}\label{eqn:fourth-kind-polynomial-iterate}
  p_k(\lambda) = \dfrac{1}{2k+1} W_k(1-2\lambda),
\end{equation}
where $W_k$ is the Chebyshev polynomial of the 4th-kind of order $k$.
Chebyshev polynomials of the 4th-kind satisfy the same recurrence as that of the first-kind, with different
initial conditions \cite{lottes2022optimal,mason1993chebyshev}:
\begin{equation}\label{eqn:fourth-kind-chebyshev}
  W_n(x) = 2xW_{n-1}(x) - W_{n-2}(x) \text{ with } W_0(x) = 1, W_1(x) = 2x+1.
\end{equation}
Following this recurrence and relaxing the $\rho(SA)=1$ assumption, a similar iterative algorithm to \cref{alg:cheby} can be derived as \cref{alg:4th-kind-cheby}.
The inclusion of the $\beta_i$ parameters in \cref{alg:4th-kind-cheby} comes as a result of optimizing
\cref{eq:gamma-defn} in the context of the \emph{multi-level} error bound from
\cref{lem:multilevel-poly-bound} \cite{lottes2022optimal}
\footnote
{
  Tabulated $\beta_i$ coefficients for the 4th-kind Chebyshev polynomials are available in \cref{tab:beta_coef}.
}
.
For the standard 4th-kind Chebyshev polynomial, $\beta_i\isdef 1$ for all $i$.
The simple smoothing iteration with $\omega=3/2$,
1st-kind Chebyshev smoothing with $\lambda_{min}=0.3\lambda_{max}$,
and the 4th-kind Chebyshev smoothing polynomials, with and without optimal $\beta_i$, are shown in \cref{fig:poly-smoothers}.
Minimizing \cref{eq:gamma-defn} with respect to $\lambda_{min}$ in the 1st-kind Chebyshev smoothing iteration is also shown.

\begin{algorithm}\small
  \caption{\small Chebyshev smoother, (Opt.) 4th-kind}
  \label{alg:4th-kind-cheby}
  \begin{algorithmic}
    \STATE{$\underline x_0 = \underline x$, $\underline r_0 = \underline b - A \underline x_0$}
    \STATE{$\underline d_0 = \dfrac{4}{3}\dfrac{1}{\lambda_{max}} S \underline r_0$}
    \FOR{$i=1,\dots,\iter-1$}
      \STATE{$  \underline x_{i} = \underline x_{i-1} + \beta_i \underline d_{i-1}, \underline r_i = \underline r_{i-1} - A \underline d_{i-1}$}
      \STATE{$  \underline d_{i} = \dfrac{2i-1}{2i+3} \underline d_{i-1} + \dfrac{8i+4}{2i+3}\dfrac{1}{\lambda_{max}}S\underline r_{i}$}
    \ENDFOR
    \STATE{$  \underline x_{\iter} =   \underline x_{\iter-1} +   \beta_k \underline d_{\iter-1}$}
    \RETURN{$  \underline x_{\iter}$}
  \end{algorithmic}
\end{algorithm}

Given $2k$ smoothing steps, what is the optimal way to distribute the number of smoothing steps
for the pre-smoother, $m$, and the post-smoother, $n$, with $m+n=2k$?
The \emph{multi-level} error bound in \cref{lem:multilevel-poly-bound}, developed by Lottes in \cite{lottes2022optimal},
provides a theoretical framework to answer the question.
\begin{lemma} \label{lem:multilevel-poly-bound}
Let the smoother iteration (on each level $j$) be given by
\[ G_j = p_{k_j}(S_jA_j) \]
where $S_j$ is SPD and scaled such that $\rho(S_jA_j) = 1$, and $p_{k_j}(x)$ is a $k_j$-order polynomial satisfying $p_{k_j}(0) = 1$ and $|p_{k_j}(x)| < 1$ for $0 < x \le 1$, possibly different on each level.
Then the V-cycle contraction factor
\begin{equation}\label{eqn:multilevel-poly-bound}
  \norm{E_{\searrow}}_{A}^2
    \leq \max_{j \in {0,\dots,\ell-1}} \dfrac{C_j}{C_j + \gamma^{-1}_j}
\end{equation}
where $C_j$ is the approximation property constant for level $j$, defined in \cref{eq:C-def},
and
\begin{equation} \label{eq:gamma-defn} \gamma_j = \sup_{0 < \lambda \le 1} \frac{\lambda \, p_{k_j}(\lambda)^2}{1 - p_{k_j}(\lambda)^2} . \end{equation}
\end{lemma}

Let us restrict the order considered in \cref{lem:multilevel-poly-bound} to the case $k_j=k$
for all $j\in{0,\dots,\ell-1}$, then $\gamma_j=\gamma$ is constant across all levels.
Further, let us define
\begin{equation}\label{eq:max-C-def}
  C \isdef \max_{j\in{0,\dots,\ell-1}} C_j.
\end{equation}
\cref{lem:multilevel-poly-bound} simplifies to
\begin{align}\label{eq:simple-multilevel-bound}
  \norm{E_{\searrow}}_{A}^2
    \leq \dfrac{C}{C + \gamma^{-1}(k)} \nonumber \\
    \isdef V(C,k).
\end{align}

Rather than directly considering the effect of the choice of $m$ and $n$ on the multigrid error contraction factor,
we instead consider the effect on the error contraction bound presented in \cref{eq:simple-multilevel-bound}.
The problem of interest, therefore, is to find $m,n$ 
\begin{align}\label{eq:optimal-mn}
  m^*, n^* & \isdef \argmin_{m,n,m+n=2k} \sqrt{V(C,m)} \cdot \sqrt{V(C,n)} \nonumber \\
  &= \argmax_{m,n,m+n=2k} C\left(\gamma^{-1}(m) + \gamma^{-1}(n)\right) + \gamma^{-1}(m)\cdot\gamma^{-1}(n),
\end{align}
where $\gamma^{-1}(k)$ is the inverse of $\gamma$, as defined in \cref{eq:gamma-defn}, for a $k$-order polynomial.

Lottes \cite{lottes2022optimal} notes that, for the simple smoother iteration with weight $\omega$,
\begin{equation}\label{eq:simple-smoother-iter-gamma}
  \gamma^{-1}_{s}(k) = 2\omega k.
\end{equation}
The 4th-kind Chebyshev polynomial, however, has
\begin{equation}\label{eq:gamma-4th-kind}
  \gamma^{-1}_{4}(k) = \dfrac{4}{3} k (k+1),
\end{equation}
while the optimized 4th-kind Chebyshev polynomial is
\begin{equation}\label{eq:gamma-opt-4th-kind}
  \gamma^{-1}_{4_{opt}}(k) = \dfrac{4}{\pi^2}(2k+1)^2 - \dfrac{2}{3}.
\end{equation}
Most notably, both 4th-kind Chebyshev polynomials types
improve the $O(k)$ simple smoother iteration to $O(k^2)$ as $k\rightarrow\infty$.
To apply \cref{lem:multilevel-poly-bound} to the 1th-kind Chebyshev polynomial,
without loss of generality, let us assume that $\lambda_{max} = 1$.
For the 1st-kind Chebyshev polynomials with fixed $\lambda_{min}$,
\begin{equation}\label{eq:gamma-1st-kind}
\gamma^{-1}_{1}(k)= \left\{
\begin{array}{ll}
      \left(T_{k}\left(\frac{\lambda_{min} + 1}{\lambda_{min} - 1}\right)\right)^{2} - 1 & k\leq k^{*}\\
      \\ %
      \frac{4 k U_{k - 1}\left(\frac{\lambda_{min} + 1}{\lambda_{min} - 1}\right)}{\left(\lambda_{min} - 1\right) T_{k}\left(\frac{\lambda_{min} + 1}{\lambda_{min} - 1}\right)} & k>k^{*} \\
\end{array} 
\right.,
\end{equation}
where $T_\xi$ and $U_\xi$ are the $\xi$th-order Chebyshev polynomials of the first and second kinds, respectively.
For $\lambda_{min}=0.1$, $k^{*}=3$.
As $k\rightarrow\infty$, \cref{eq:gamma-1st-kind} scales as
\begin{equation}\label{eq:gamma-1st-kind-asm}
  \gamma^{-1}_{1}(k) \sim \sqrt{\dfrac{1}{\lambda_{min}}} k.
\end{equation}
However, by choosing $\lambda_{min}$ such that $\gamma^{-1}$ is maximized,
\begin{equation}\label{eq:gamma-opt-1st-kind}
  \gamma^{-1}_{1_{opt}}(k) \sim 2.38k^{1.73}
\end{equation}
for large $k$.
As an aside, a correlation for $\lambda_{min}^{*}$ with 1\% relative error
and 0.1\% absolute error for $k\in[1,50]$ is given by
\begin{equation}\label{eq:opt-lambda-min}
  \lambda_{min}^* \approx \dfrac{1.69}{k^{1.68} + 2.11k + 1.98}.
\end{equation}

Due to symmetry of \cref{eq:optimal-mn}, the error bound for $(m,n)$ is the same as that of $(n,m)$.
To start to assess \cref{eq:optimal-mn}, we first consider the case with $m=n=k$ compared to $m=2k$, $n=0$.
For the simple smoother with fixed $\omega$, the error-bound is minimized with $m=n$.
For the 4th-kind Chebyshev polynomial, $m=2k$, $n=0$ outperforms the symmetric $m=n=k$ if and only if
\begin{equation}\label{eq:4th-kind-condition-C}
  C > \frac{2 \left(k + 1\right)^{2}}{3}
\end{equation}
The optimized 4th-kind Chebyshev polynomial has a similar condition, provided
\begin{equation}\label{eq:opt-4th-kind-condition-C}
  C > \frac{2 \left(6 \left(2 k + 1\right)^{2} - \pi^{2}\right)^{2}}{3 \pi^{2} \left(- 12 \left(2 k + 1\right)^{2} + 6 \left(4 k + 1\right)^{2} + \pi^{2}\right)}
\end{equation}
For the 1st-kind Chebyshev polynomial with fixed $\lambda_{min}$,
\begin{equation}\label{eq:1st-kind-condition-C}
  C \gtrapprox 1.55 e^{1.45k}
\end{equation}
as $k\rightarrow\infty$.
However, for $k\in[1,3]$, $m=2k$, $n=0$ outperforms $m=n=k$, irrespective of $C$.
Lastly, for the 1st-kind Chebyshev polynomial with $\lambda_{min}$ chosen to maximize $\gamma^{-1}$,
\begin{equation}\label{eq:opt-1st-kind-condition-C}
  C \gtrapprox 1.81 k^{1.73}.
\end{equation}

With the regions where $m=2k$, $n=0$ outperforms $m=n=k$ determined,
the authors wish to solve the more general optimization problem in \cref{eq:optimal-mn}.
In lieu of a general solution, however, the authors opt to prove that the optimal $m$ and $n$
occurs at either $m=2k$, $n=0$ or $m=n=k$ for all $C>0$.
To do so, the authors utilized sympy \cite{meurer2017sympy} to solve for the region
where a tentative $(\tilde m,\tilde n)$ outperforms $m=2k$, $n=0$ and the region where $(\tilde m,\tilde n)$ outperforms $m=n=k$.
The intersection of these two regions, therefore, is the region where $(\tilde m,\tilde n)$ outperforms
both $m=2k$, $n=0$ and $m=n=k$.
With exception to $(4,2)$ and $(5,1)$ for $C<4$ in the case of 1st-kind Chebyshev polynomial smoothing with $\lambda_{min}=0.1$,
the authors found that the optimal $(m,n)$ is either $(m=2k,n=0)$ or $(m=n=k)$ for all $C>0$ for the smoothers considered,
up to $k=50$.

A natural next question is on the expected improvement using the $(2k,0)$ scheme
over the symmetric $(k,k)$ scheme.
This is done by taking the ratio of the error bound \cref{eq:optimal-mn} with $(k,k)$ and $(2k,0)$
for the various polynomial smoothers considered.
For the 1st-kind Chebyshev polynomial with $\lambda_{min}=0.1$, the expected improvement in the multigrid
convergence is no more than 6\% for $k=1$, and quickly decreases for larger $k$.
Optimizing the bound with respect to $\lambda_{min}$, however, the 1st-kind Chebyshev polynomial
with $(2k,0)$ can outperform the symmetric $(k,k)$ scheme by 9\%.
For both the 4th-kind and optimized 4th-kind Chebyshev polynomials, however,
the $(2k,0)$ scheme outperforms the symmetric $(k,k)$ scheme by 15\% as $k\rightarrow\infty$.
A numerical example demonstrating the applicability of this analysis based on \cref{lem:multilevel-poly-bound}
and \cref{eq:optimal-mn} is given in \cref{sec:fd-gmg}.

\section{Finite Differences with Geometric Multigrid Poisson}
\label{sec:fd-gmg}
\def\dO{\partial \Omega}
\def\Oh{{\hat \Omega}}
\def\uu{{\underline u}}
\def\ur{{\underline r}}
\def\us{{\underline s}}
\def\uf{{\underline f}}
\def\tlam {{\tilde \lambda}}
\def\interp {{P}}

The Poisson equation \cref{eqn:poisson} is considered with $d=2$.
Let $\Omega:=[0,L_x] \times [0,L_y]$ be the domain of interest, with the boundary condition $u|_{\dO} = 0$.
A finite difference grid of $(n+1)\times (n+1)$ points is considered, $n=128$.
For the purposes of this study, $u(x,y) = \sin\left({3\pi x}/{L_x}\right)\sin\left({4\pi y}/{L_y}\right)+g$,
where $g$ is the same random vector with $g|_{\dO}=0$.
$L_x\ge 1$ is varied, while $L_y$ is fixed at unity.
A geometric multigrid V-cycle with Chebyshev-accelerated Jacobi smoothing is employed as a preconditioner for GMRES(20).
For the comparison with the bounds presented in \cref{lem:multilevel-poly-bound},
multigrid is also considered as a solver.
Each coarser level is discretized as $(n_{c} + 1)\times (n_{c} + 1)$,
with $n_{c} = n/2$, as well as aggressive coarsening with $n_{c} = n/8$.
This is repeated until there is only a single degree of freedom.
As this case is meant to represent our target problem for \emph{unstructured}
multigrid, semi-coarsening is not employed.
On the coarsest level, the single degree of freedom system is solved exactly.
We choose a relative residual reduction of $10^{-6}$ as the stopping criterion.
A variety of smoothing orders are considered for all orders $m$ pre-smoothing
and orders $n$ post-smoothing, with $m+n=2k$ up to $k=10$.

\begin{figure}[h]
  \centering
  \def\svgwidth{\textwidth}
  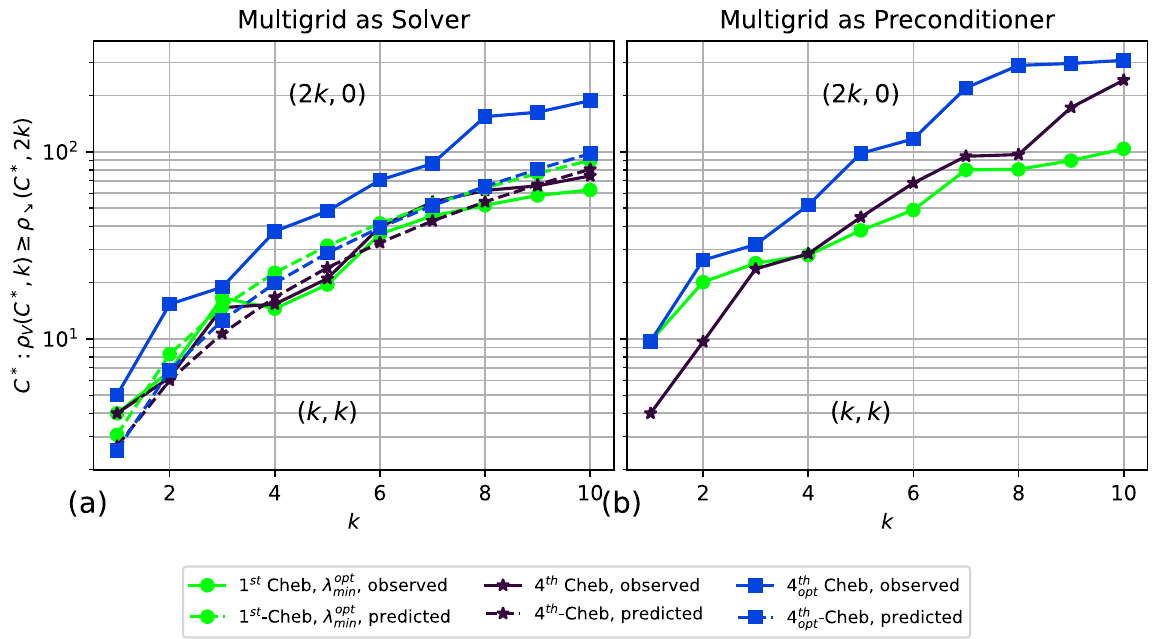
  \caption{\small 
    \label{fig:critical-c}
    Critical $C^*$ at which $(2k,0)$ converges faster than $(k,k)$ for $C > C^{*}$.
  }
\end{figure}

\begin{figure}[h]
  \centering

  \def\svgwidth{\textwidth}
  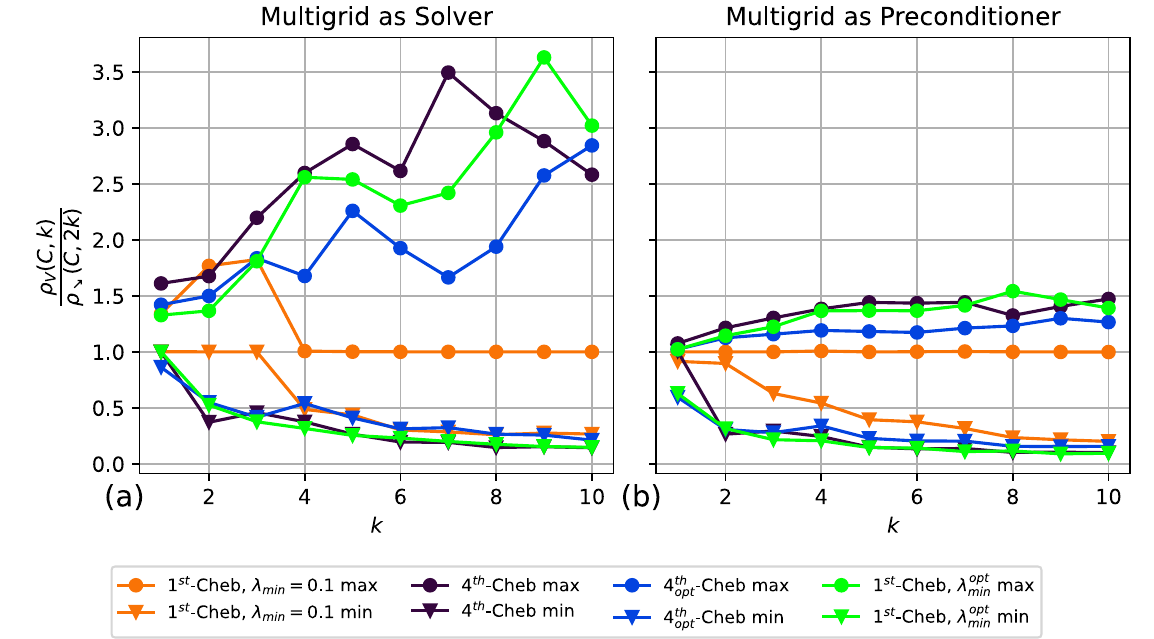
  \caption{\small 
    \label{fig:min-max-ratios}
    Max (min) error contraction rate ratios for $(2k,0)$ and $(k,k)$ smoothing schemes.
  }
\end{figure}

Convergence results are based on the observed average solver convergence rate
\begin{equation}\label{eq:avg-convergence-rate}
  \rho = \exp\left({\dfrac{1}{N} \log \dfrac{\norm{\underline r_N}}{\norm{\underline r_0}}}\right),
\end{equation}
where $N$ is the number of iterations.
The first question, as posed in the optimization problem \cref{eq:optimal-mn},
is regarding the optimal choice of V-cycle smoothing $(m,n)$ with $m+n=2k$.
Recall, however, that the authors demonstrated the optimality of
either the one-sided $(2k,0)$ or symmetric $(k,k)$ V-cycles for all orders
up to $k=50$, with a few exceptions for small values of $C$.
The authors observe that, as predicted in \cref{eq:simple-multilevel-bound,eq:optimal-mn},
the convergence rate of $(m,n)$ is nearly equivalent to $(n,m)$.
Further, with few exceptions, either the one-sided $(2k,0)$ or symmetric $(k,k)$ V-cycles
yielded the smallest convergence rate, as defined in \cref{eq:avg-convergence-rate}.

With our consideration now on the use of the symmetric $(k,k)$ or one-sided $(2k,0)$ V-cycles,
we can now confirm our theoretical prediction of when to apply these methods,
as outlined in \cref{eq:4th-kind-condition-C,eq:opt-4th-kind-condition-C,eq:1st-kind-condition-C,eq:opt-1st-kind-condition-C}.
The convergence rate for the symmetric $(k,k)$ V-cycle is denoted $\rho_V(C,k)$, while the one-sided $(2k,0)$ V-cycle is $\rho_{\searrow}(C,2k)$.
Values of $(C,k)$ at which $\rho_V(C,k) \geq \rho_{\searrow}(C,2k)$ are shown in \cref{fig:critical-c}.
Predicted results are from \cref{eq:4th-kind-condition-C,eq:opt-4th-kind-condition-C,eq:1st-kind-condition-C,eq:opt-1st-kind-condition-C},
which are only applicable to multigrid as a \emph{solver}, not \emph{preconditioner}.
Observed results are from 2D finite difference example.
At a given $k$, $C > C^{*}$ indicates that the one-sided $(2k,0)$ V-cycle yields a lower error bound than the symmetric $(k,k)$ V-cycle.
Conversely, $C < C^{*}$, indicates that the symmetric V-cycle yields better convergence.
For the 1st-kind with optimized $\lambda_{min}$ coefficient,
4th-kind, and optimized 4th-kind Chebyshev smoothers,
the predicted domain in which to apply the one-sided $(2k,0)$ V-cycle over
the symmetric $(k,k)$ V-cycle shows agreement
with the results obtained by experiment using multigrid as a solver (\cref{fig:critical-c}a).
As noted by \cref{eq:1st-kind-condition-C}, when $k>3$, the 1st-kind Chebyshev smoother
does not benefit from applying the one-sided $(2k,0)$ V-cycle over the symmetric $(k,k)$ V-cycle.
At the same time, however, for $k\leq 3$, the one-sided $(2k,0)$ V-cycle
outperforms the symmetric $(k,k)$ V-cycle, irrespective of $C$.
Despite the applicability of \cref{lem:multilevel-poly-bound} being limited to multigrid as a solver,
the predicted domain in which to apply the one-sided $(2k,0)$ V-cycle \emph{as a preconditioner} is similar to that of using multigrid \emph{as a solver} (\cref{fig:critical-c}b).

The maximum and minimum ratio of the error contraction rates, along with
the predicted performance, are shown in \cref{fig:min-max-ratios}.
While the predicted performance benefit of applying the one-sided V-cycle approach is limited,
nevertheless, \cref{fig:min-max-ratios}a and \cref{fig:min-max-ratios}b demonstrate
that the one-sided V-cycle offers an improvement compared to the symmetric V-cycle
for problems with moderate values of $C$.
However, when applied to relatively easy problems ($C\approx 1$), the one-sided V-cycle
is a poor choice.

\begin{figure}
  \centering
  \def\svgwidth{\textwidth}
  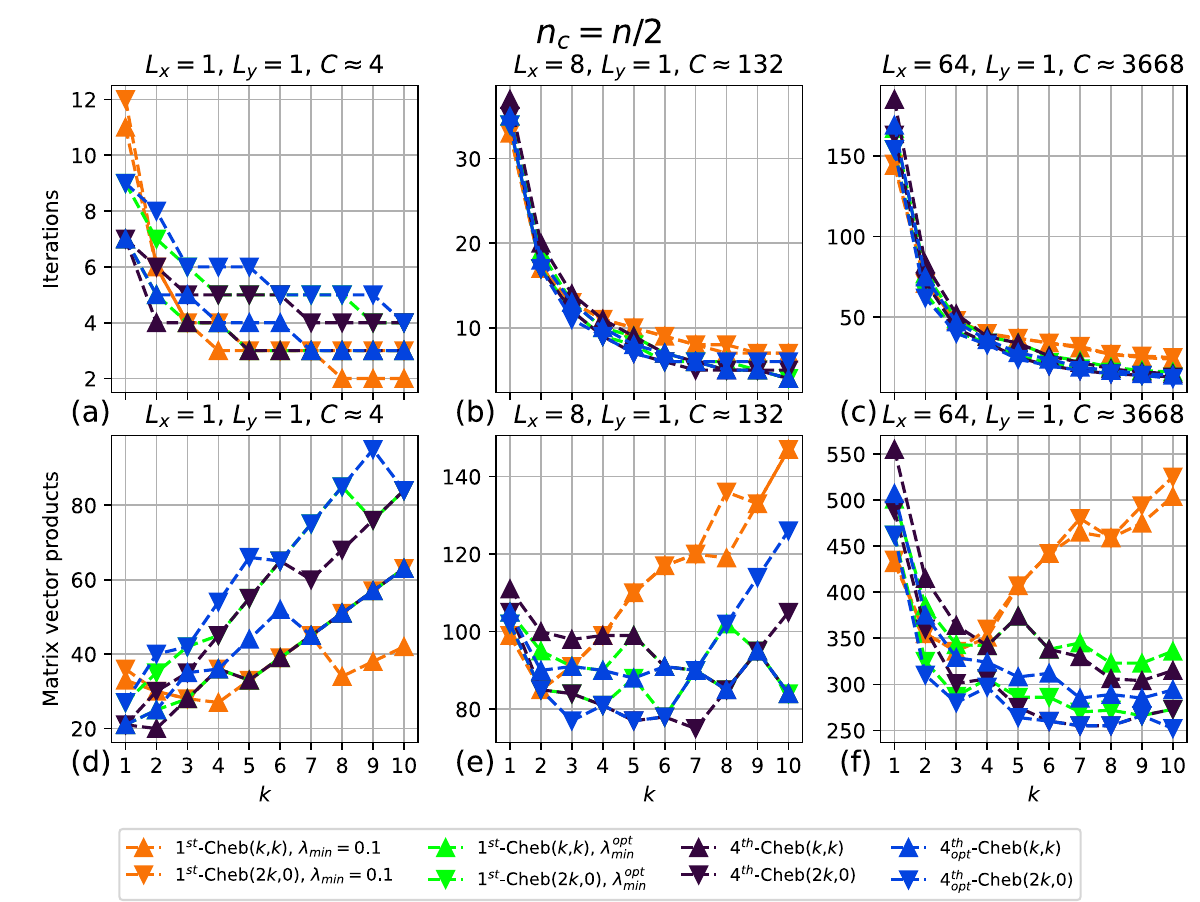
  \caption{\small \label{fig:fd-factor-2} FD, $n_{c} = n/2$. Multigrid as preconditioner for GMRES(20).}
\end{figure}

\begin{figure}
  \centering
  \def\svgwidth{\textwidth}
  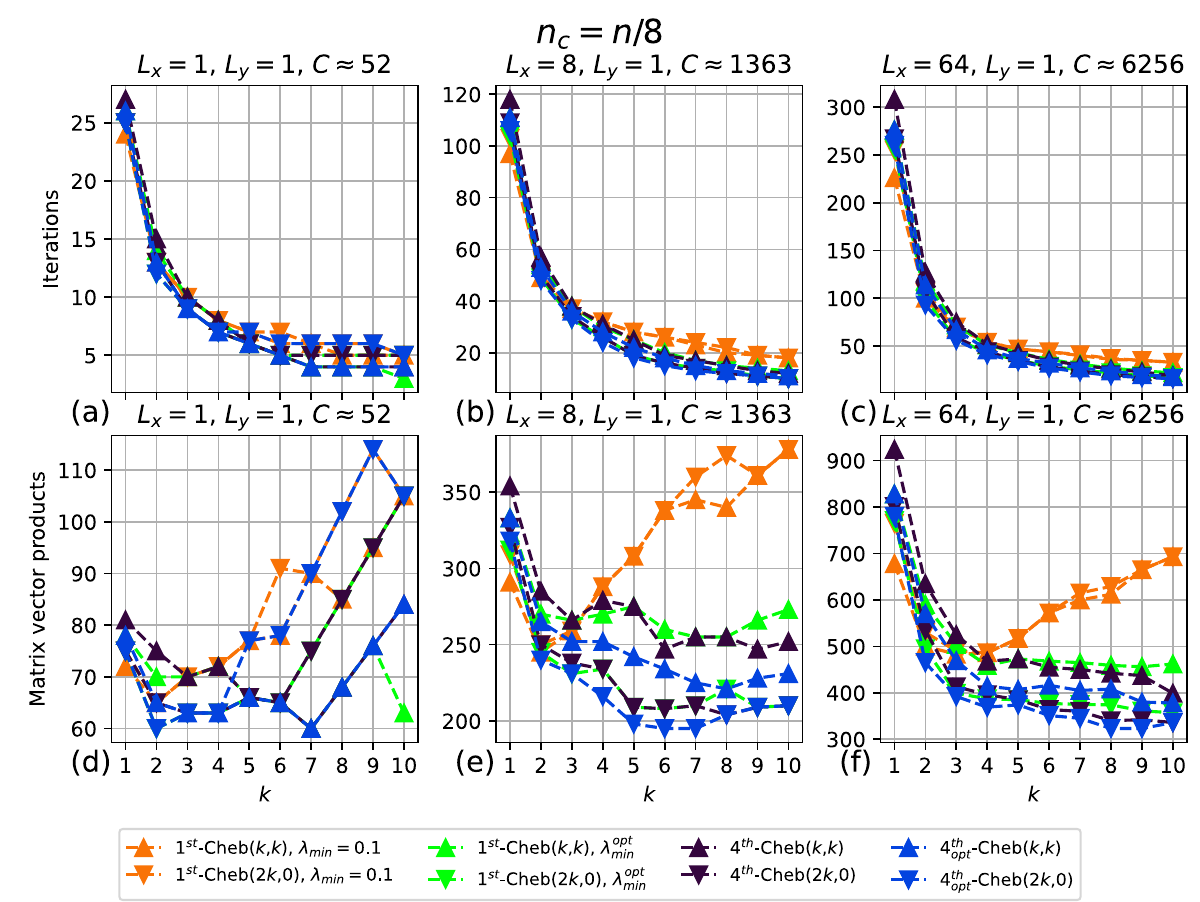
  \caption{\small \label{fig:fd-factor-8} FD, $n_{c} = n/8$. Multigrid as preconditioner for GMRES(20).}
\end{figure}

Results for $n_{c} = n/2$ are shown in \cref{fig:fd-factor-2}, $n_c = n/8$ in \cref{fig:fd-factor-8}.
In both, multigrid is used to precondition the GMRES(20) solver.
In the case $n_c = n/2$ and $L_x=1$ ($C\approx 4$), shown in \cref{fig:fd-factor-2}a,d,
the work required, as measured by fine-grid matrix vector products,
is minimized for relatively low orders.
Further, applying the bounds shown in \cref{fig:critical-c}a,
this is the scenario in which the one-sided V-cycle approach \emph{should not}
be utilized.
However, even at moderate grid aspect ratios, such as $L_x=8$, $C\approx 132$ becomes
large enough to justify the usage of the one-sided V-cycle approach, especially for
the 4th-kind Chebyshev smoothers (\cref{fig:fd-factor-2}b,e) at moderate orders
$k=6, \tilde k = 12$.
This effect becomes even more apparent when $L_x=64$ ($C\approx 3668$), shown in \cref{fig:fd-factor-2}c,f.
In the results for $n_c=n/2$, the 1st-kind with optimized $\lambda_{min}$, 4th-kind, and optimized 4th-kind Chebyshev smoothers
greatly outperform the 1st-kind Chebyshev smoother, \emph{especially} at high orders.
Further, it is observed that both the iteration count and \emph{total work} in matrix-vector products
is minimized at high orders for the 1st-kind with optimized $\lambda_{min}$,
4th-kind, and optimized 4th-kind Chebyshev smoothers.
This is not the case for the 1st-kind Chebyshev smoother, however.
This effect of lowering both the iteration count and work has
the additional benefit of reducing the number
of coarse grid solves required for each iteration, whose cost is not factored in this analysis.

In the aggressive coarsening case ($n_c = n/8$), both the number of matrix-vector products and iterations are reduced
through using higher-order Chebyshev smooothing irrespective of the grid aspect ratio
for all but the standard 1st-kind Chebyshev smoothers.
Secondly, the one-sided $(2k,0)$ V-cycle approaches generally outperforms the full, symmetric $(k,k)$ V-cycle approach.
There are two important implications from this result:
first, there exists considerable benefit in transitioning
a multigrid solver to either construct high-quality estimates for $\lambda_{min}^{opt}$ for the 1st-kind Chebyshev smoother,
such as those provided in correlation \cref{eq:opt-lambda-min},
or utilize one of the 4th-kind Chebyshev smoothers;
and second, additional performance can be achieved by using the one-sided $(2k,0)$ V-cycle approach at the expense of symmetry,
the implication of which is further discussed in \cref{sec:loss-of-symmetry}.

Let us also consider the trade-offs associated with using $n_c = n/2$ and $n_c = n/8$ as coarsening strategies.
Despite improvements in the convergence rate from using a different V-cycle approach or Chebyshev smoother,
the number of matrix-vector products required by the solve is greater for the aggressive coarsening case.
However, the grid complexity
\begin{equation}\label{eq:mg-cycle-complexity}
  \dfrac{\sum_{j=0}^{\ell}\nnz{A_j}}{\nnz{A_0}}
\end{equation}
for $n_c = n/2$ is 1.568, while for $n_c = n/8$ it is only 1.023.
Therefore, the \emph{amount} of work done per cycle is less for the aggressive coarsening case.
An additional benefit of the aggressive coarsening strategy is decreasing the number of multigrid levels
from 7 with $n_c=n/2$ to 3 with $n_c=n/8$.
This feature is especially attractive in a parallel computing context
where each additional multigrid level requires additional communication cost.
Significant effort has been spent on improving the parallel scalability of multigrid methods,
especially in the AMG context, see \cite{stuben1999algebraic,falgout2014non,bienz2020reducing}.
One strategy to achieve better scalability is to rely on aggressive coarsening strategies,
which require more robust smoothers, such as the Chebyshev smoothers discussed here.

The results shown in \cref{fig:fd-factor-2} and \cref{fig:fd-factor-8} are summarized
in \cref{tab:fastest_solver_fd}.
The solver configuration yielding the lowest number of matrix-vector products is listed for each case.
$4^{th}_{opt}$-Cheb, Jacobi(20,0), for example, denotes a 20th order Chebyshev smoother
of the optimized 4th-kind, using one-sided smoothing (thereby, having the same cost per iteration as order 9
with the symmetric V-cycle).
This solver configuration yields the lowest number of matrix-vector products
for $L_x=128$ with aggressive coarsening.
We see that, for high aspect ratio grids and aggressive coarsening,
the one-sided $(2k,0)$ V-cycle offers superior performance
to the symmetric $(k,k)$ V-cycle approach.

\begin{table}
\centering
\caption{\label{tab:fastest_solver_fd}Solver configuration with lowest number of matrix-vector products for finite different geometric multigrid}
\begin{tabular}{llllrl}
\toprule
$L_x$ & Coarsening & Complexity &                                   Solver &  Mat-Vec & Iterations \\
\midrule
     1 &      2 & 1.568 & $4^{th}$-Cheb, Jacobi(2, 2) &     20 &          4 \\
     8 &      2 & 1.568 & $4^{th}$-Cheb, Jacobi(14, 0) &     75 &          5 \\
    64 &      2 & 1.568 & $4^{th}_{opt}$-Cheb, Jacobi(20, 0) &    252 &         12 \\
   128 &      2 & 1.568 & $4^{th}_{opt}$-Cheb, Jacobi(20, 0) &    252 &         12 \\
\midrule
     1 &      8 & 1.023 &      $4^{th}$-Cheb, Jacobi(7, 7) &     60 &          4 \\
     8 &      8 & 1.023 & $4^{th}_{opt}$-Cheb, Jacobi(14, 0) &    195 &         13 \\
    64 &      8 & 1.023 & $4^{th}_{opt}$-Cheb, Jacobi(18, 0) &    323 &         17 \\
   128 &      8 & 1.023 & $4^{th}_{opt}$-Cheb, Jacobi(20, 0) &    294 &         14 \\
\bottomrule
\end{tabular}
\end{table}

\section{High Order Preconditioners}
\label{sec:pmg-preconditioners}
Let us now consider preconditioners for the Poisson equation \cref{eqn:poisson} for $d=3$
arising from the spectral element method (SEM) discretization of the incompressible NS equation,
which typically encompasses the majority of the solution time in nekRS \cite{fischer2022nekrs}.
Two classes of preconditioners prove most effective: $p$-geometric multigrid (pMG)
and low-order discretizations.

\subsection{$p$-Geometric Multigrid} \label{sec:pmg}

$p$-geometric multigrid (pMG) is used as a preconditioner for the Pressure poisson equation
discretized using the spectral element method.
Since the multigrid hierarchy is constructed by varying the polynomial order, $p$, of each level,
geometric coarsening can be robustly applied up to $p=1$, even in the \emph{unstructured} case.
Typical multigrid schedules for $p=7$, for example, are V-cycles with orders $(7,5,3,1)$ or $(7,3,1)$.
In addition to the Chebyshev-accelerated Jacobi smoother,
Chebyshev-accelerated Schwarz smoothers are also considered, as discussed in \cite{phillips2022tuning}.
The Schwarz smoothers are outlined below.

The SE-based additive Schwarz method (ASM) presented in
\cite{lottes2005hybrid,loisel2008hybrid} solves local Poisson
problems on subdomains that are extensions of the spectral elements.
The formal definition of the ASM preconditioner (or, in this case, pMG smoother) is
\begin{equation} \label{eq:asm}
  S_{ASM} \underline r = \sum_{e=1}^E W_e R_e^T {\bar A}_e^{-1} R_e \underline r,
\end{equation}
where $R_e$ is the restriction matrix that extracts nodal values
of the residual vector that correspond to each overlapping domain.
To improve the smoothing properties of the ASM, we introduce
the diagonal weight matrix, $W_e$, which scales each nodal value by
the inverse of the number of subdomains that share that node.
Although it compromises symmetry, post-multiplication by $W_e$ was found to yield
superior results to pre- and post-multiplication by $W_e^{\frac{1}{2}}$
\cite{stiller2017nonuniformly,lottes2005hybrid}.

In a standard Galerkin ASM formulation, one would use ${\bar A}_{e} = R_e A
R_e^T$, but such an approach would compromise the $O(p^3)$ storage complexity
of the SE method. To construct fast inverses for ${\bar A}_e$, we approximate
each deformed element as a simple box-like geometry,
These boxes are then extended by a single 
degree-of-freedom in each spatial dimension to form overlapping subdomains with
${\bar p}^3=(p+3)^3$ interior degrees-of-freedom in each domain.
The approximate box domain enables the use of the fast diagonalization method
(FDM) to solve for each of the overlapping subdomains, which can be applied in
$O(Ep^4)$ time in $\RR^3$. The
extended-box Poisson operator is
\begin{equation*}\label{eq:tensor-prod-poisson}
{\bar A} = 
B_z \otimes B_y \otimes A_x+B_z \otimes A_y \otimes B_x+A_z \otimes B_y \otimes B_x,
\end{equation*}
where each $B_*,A_*$ represents the extended 1D mass-stiffness matrix pairs
along the given dimension \cite{dfm02}. The FDM begins with a
preprocessing step of solving a series of small,
${\bar p}\times {\bar p}$, generalized eigenvalue problems, 
\begin{equation*}\label{eq:generalized-eig}
A_* \underline{s}_i = \lambda_i B_* \underline{s}_i
\end{equation*}
and defining
$S_*=(\underline{s}_1\ldots\underline{s}_{\bar p})$ and
$\Lambda_*=\operatorname{diag}(\lambda_i)$, to yield the
similarity transforms
\begin{equation*}\label{eq:similarity}
S_*^T A_* S_* = \Lambda_*, \;\; S_*^T B_* S_* = I.
\end{equation*}
From these, the inverse of the local Schwarz operator is \\[-2.5ex]
\begin{equation*}\label{eq:inverse-tensor-prod-poisson}
{\bar A}^{-1} = 
(S_z\otimes S_y\otimes S_x) D^{-1} (S_z^T \otimes S_y^T \otimes S_x^T),
\end{equation*}
where $D$ is a diagonal matrix defined as
\begin{equation*}
D\isdef I\otimes I\otimes \Lambda_x+I\otimes \Lambda_y \otimes I+\Lambda_z \otimes I \otimes I.
\end{equation*}
This process is repeated for each element, at each 
multigrid level save for the coarsest one. Note that the
per-element storage is only $3{\bar p}^2$ for the $S_*$ matrices and ${\bar
p}^3$ for $D$. At each multigrid level, the local subdomain solves are used as
a smoother. On the coarsest level ($p=1$), however, BoomerAMG
\cite{yang2002boomeramg} is used to solve the system on the CPU with the same
parameters described in \cref{sec:semfem}, except using
Chebyshev smoothing.
A single BoomerAMG V-cycle iteration is used in the coarse-grid solve.

Presently, we also consider a restrictive additive Schwarz (RAS) version
of (\ref{eq:asm}), wherein overlapping values are not added after the action 
of the local FDM solve, following \cite{cai1999restricted}.
RAS has the added benefit of reducing the amount of communication required in the smoother.
Similar to ASM, RAS is \emph{non-symmetric}.
Attempts to symmetrize the operator tend to have a negative impact on the convergence rate \cite{cai1999restricted}.
The smoother $S_{ASM}$ or $S_{RAS}$ described here is then used as the smoother
considered in Chebyshev-acceleration, as shown in \cref{alg:4th-kind-cheby} and
\cref{alg:cheby}.
The multigrid schedule for the cases in \cref{sec:cases} are $(7,5,3,1)$
and $(7,3,1)$ for Chebyshev-Jacobi and Chebyshev-Schwarz respectively.

While the weighted ASM method in \cref{eq:asm} and RAS method
are not SPD, experimental evidence suggests that the complex spectra
of $S_{ASM}A$ and $S_{RAS}A$ have small imaginary part.
Following Manteuffel \cite{manteuffel1977tchebychev},
\cref{alg:cheby} can be adapted to the non-SPD case by
replacing the interval $[\lambda_{min},\lambda_{max}]$ with
an ellipse centered about $\theta$ with a focal distance $\delta$
that encloses the convex hull of the spectra.
However, since the eigenvalues that have an imaginary part
have moduli much less than $\lambda_{max}(SA)$,
adapting the algorithm to the non-SPD case is not necessary.
Numerous numerical results (\cite{fischer2022nekrs,min2022optimization,reger2022toward,phillips2022tuning}, etc.)
confirm the efficacy of this approach,
as pMG with Chebyshev-accelerated ASM smoothing is the default preconditioner in nekRS.

\subsection{Preconditioning with Low Order Operators} \label{sec:semfem}
In \cite{orszag1979spectral}, Orszag suggested that constructing a sparse
preconditioner based on the low-order discretizations with nodes coinciding
with those of the high-order discretization would yield bounded condition
numbers and, under certain constraints, can yield $\kappa(M^{-1}A) \sim \pi^2 /
4$ for second-order Dirichlet problems. This observation has led to the
development of preconditioning techniques based on solving the resulting low-order system
\cite{pazner2020efficient,olson2007algebraic,bello2019scalable,canuto2010finite}.

In the current work, we employ the same low-order discretization
considered in \cite{bello2019scalable}. Each of the vertices of
the hexahedral element is used to form one low-order, tetrahedral element,
resulting in a total of eight low-order elements for each GLL sub-volume
in each of the high-order hexahedral elements.
This low-order discretization is then used to form the sparse operator, $A_F$.
The so-called weak preconditioner, $A_F^{-1}$, is used to precondition the system.
boomerAMG \cite{yang2002boomeramg},
is used with the following setup to solve the low order system: \\[-2.5ex]
\begin{itemize}
  \item PMIS coarsening                              \\[-2.5ex]
  \item 0.25 strength threshold                      \\[-2.5ex]
  \item Extended + i interpolation ($p_{max}=4$)     \\[-2.5ex]
  \item L1 Jacobi smoothing                          \\[-2.5ex]
  \item One V-cycle for preconditioning              \\[-2.5ex]
  \item Smoothing on the coarsest level              \\[-2.5ex]
\end{itemize}
We denote this preconditioning strategy as SEMFEM.

\section{High Order Cases}
\label{sec:cases}
\def\dt{ \Delta t }
We describe four model problems that are used to test the high-order $p$-multigrid (pMG) preconditioners.
The first is a stand-alone Poisson solve with the Kershaw mesh ($\varepsilon=1,0.3,0.05$).
The others are modest-scale NS problems, where the pressure Poisson
problem is solved over 2,000 timesteps.
Details regarding the spectral element discretization are given in \cite{phillips2022tuning}.
The problem sizes are listed in
\cref{table:problem-sizes} and range from small ($n$=16M points) to moderate
($n$=180M).\footnote{Larger cases for full-scale runs on Summit with 
$n$=51B are reported in \cite{min2022optimization}.}

\subsection{Poisson} \label{sec:kershaw}

\begin{figure}
{\setlength{\unitlength}{\textwidth}
\begin{picture}(1,0.33)
   \put(0.00,0.025){\includegraphics[width=0.33\textwidth, height=0.5\unitlength, keepaspectratio]{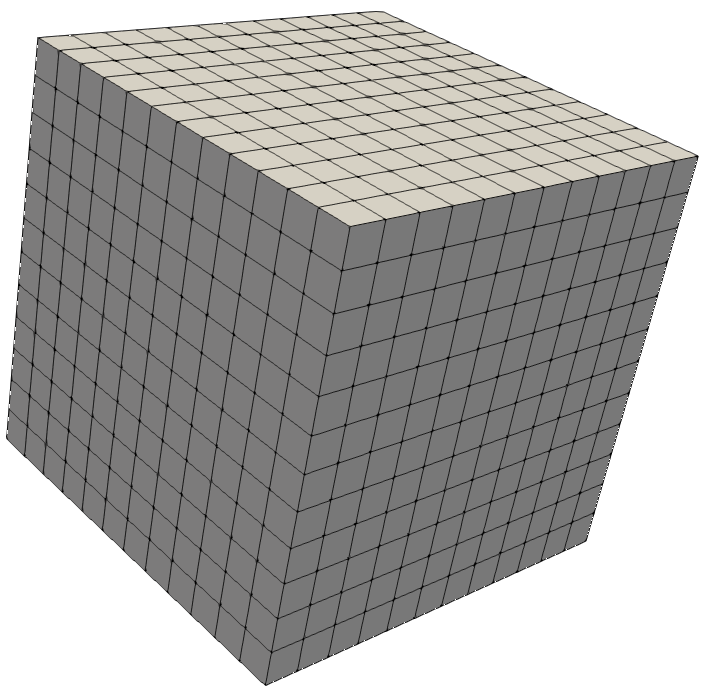}}
   \put(0.33,0.025){\includegraphics[width=0.33\textwidth, height=0.5\unitlength, keepaspectratio]{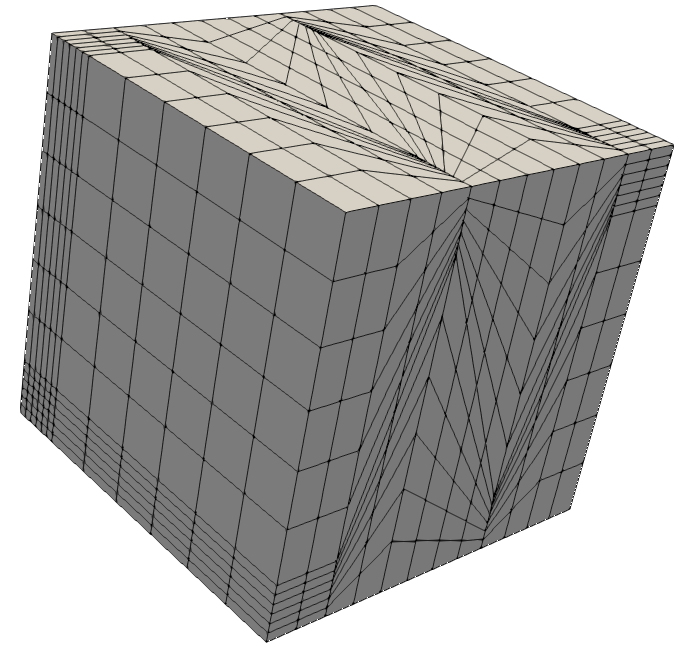}}
   \put(0.67,0.025){\includegraphics[width=0.33\textwidth, height=0.5\unitlength, keepaspectratio]{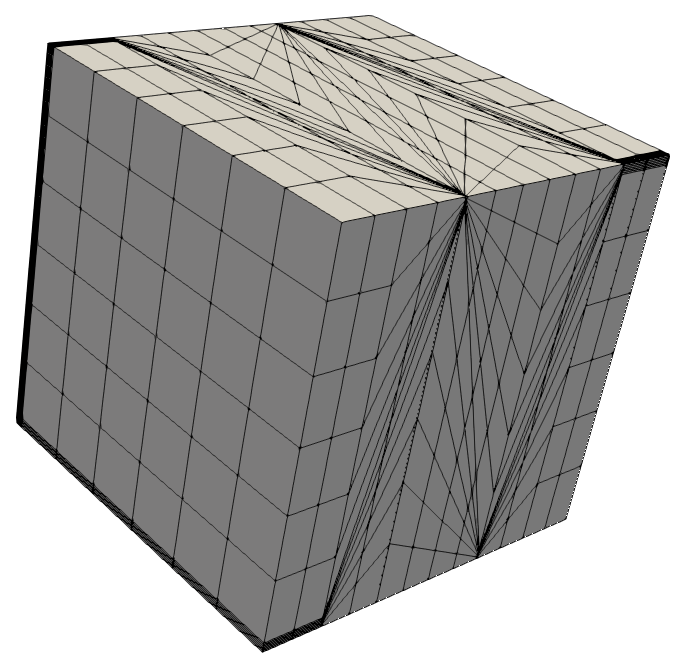}}

   \put(-0.05,0.05){$\varepsilon = 1.0$}
   \put(0.25,0.05){$\varepsilon = 0.3$}
   \put(0.60,0.05){$\varepsilon = 0.05$}

\end{picture}}
\caption{\small 
  \label{fig:kershaws}
  Kershaw, $E=12^3,p=1$.
}
\end{figure}
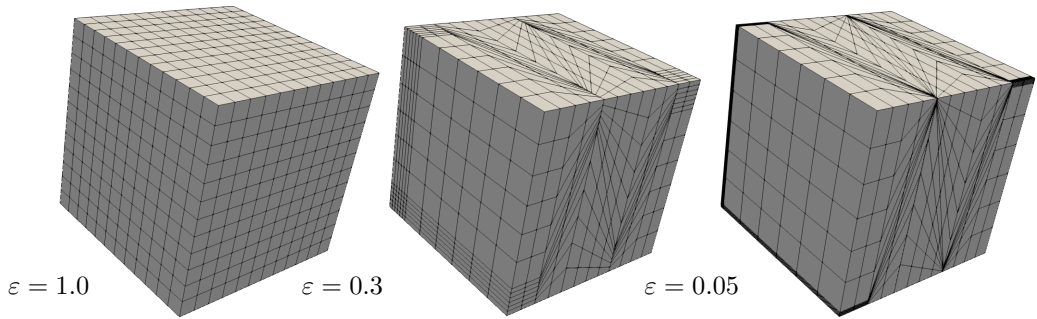

The Kershaw family of meshes \cite{kolevceed,kershaw1981differencing}
has been proposed as the basis for a high-order Poisson-solver benchmark by
Center for Efficient Exascale Discretization (CEED) within the DOE Exascale
Computing Project (ECP). This family is parameterized by an anisotropy measure,
$\varepsilon =\varepsilon_y=\varepsilon_z \in (0,1]$, that determines the
degree of deformation in the $y$ and $z$ directions. As $\varepsilon$
decreases, the mesh deformation and aspect ratio increase along with it. 
The Kershaw mesh is shown in \cref{fig:kershaws} for
$\varepsilon=1,0.3,0.05$.
The domain $\Omega = \left[-1/2, 1/2\right]^3$ with Dirichlet boundary conditions on $\partial\Omega$.
The right hand side for \cref{eqn:poisson} is set to
\begin{equation}\label{eq:kershaw-rhs}
  f(x,y,z) = 3\pi^2 \sin{(\pi x)}\sin{(\pi y)}\sin{(\pi z)}.
\end{equation}
The linear solver terminates after reaching a relative residual reduction of
$10^{-8}$.
For the $1^{st}$-Cheb, $\lambda_{min}^{opt}$ smoother method,
a random right-hand side is used to tune $\lambda_{min}$,
rather than using the correlation from \cref{eq:opt-lambda-min}.
The solver used is GMRES(30) preconditioned with a single pMG V-cycle.
Since this test case solves the Poisson equation, there is no
timestepper needed for the model problem.

\subsection{Navier-Stokes} \label{sec:navier-stokes}
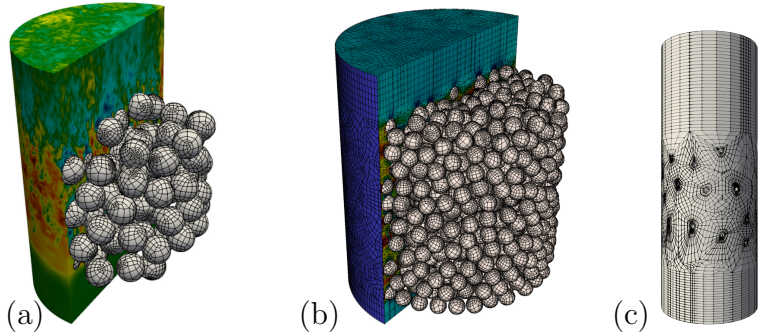
\begin{figure}
  \centering
  {\setlength{\unitlength}{\textwidth}
    \begin{picture}(1,0.33)(0,0)
      \put(0.00,0.0){\includegraphics[width=0.3\textwidth, height=0.33\unitlength, keepaspectratio]{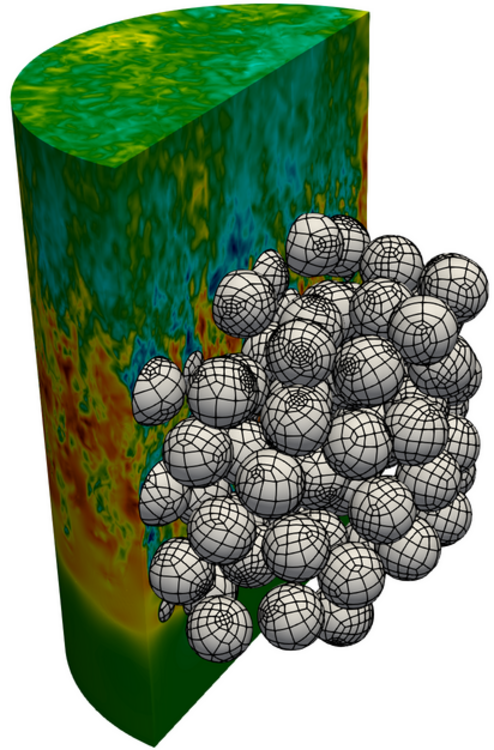}}
      \put(0.33,0.0){\includegraphics[width=0.3\textwidth, height=0.33\unitlength, keepaspectratio]{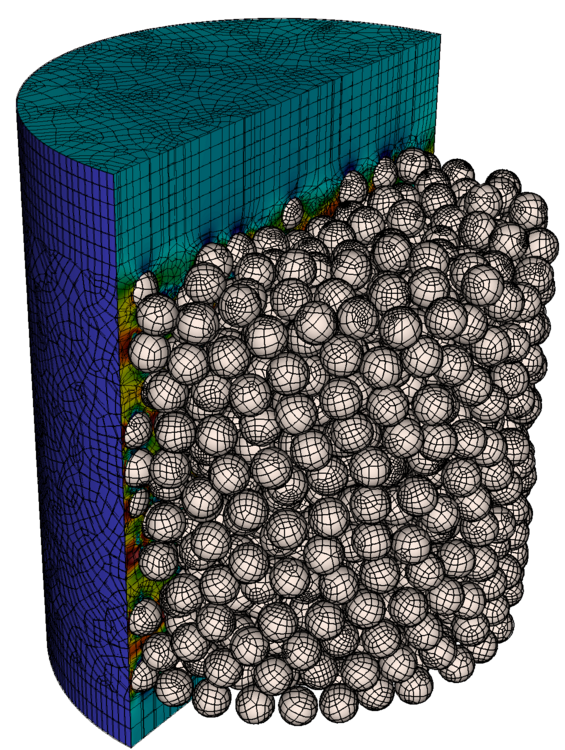}}
      \put(0.67,0.0){\includegraphics[width=0.3\textwidth, height=0.33\unitlength, keepaspectratio, angle=90]{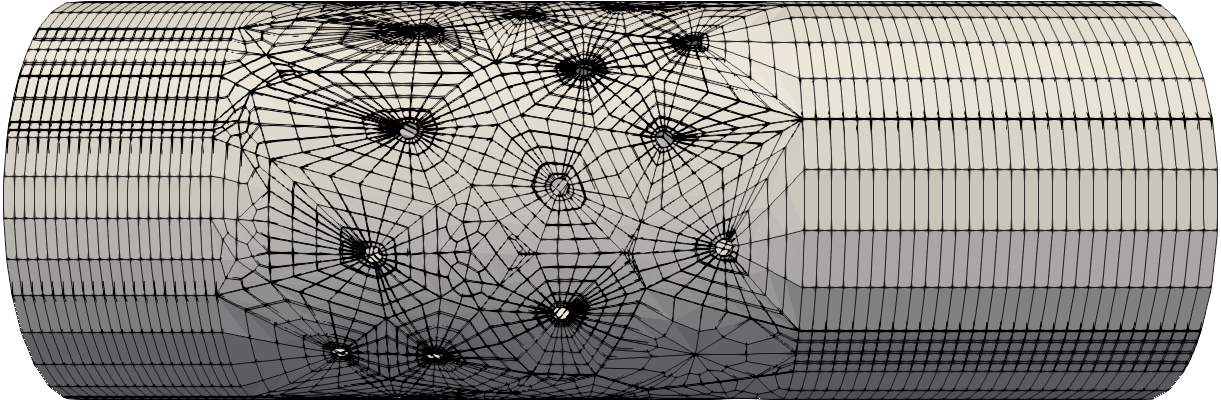}}

      \put(0.00,0.0){\large (a)}
      \put(0.30,0.0){\large (b)}
      \put(0.62,0.0){\large (c)}

    \end{picture}}
  \caption{\small 
    Navier-Stokes cases: pebble-beds with 
    (a) 146, (b) 1568, and (c) 67 spheres
    \label{fig:ns_cases}
  }
\end{figure}

\begin{table}
\centering
\caption{
  Problem discretization parameters.
  \label{table:problem-sizes}
}
\begin{tabular}{lrrrrr}
  \toprule
  Case & $E$ & $p$ & $n$ & $P$ & $\dfrac{n}{P}$ \\
  \midrule
  Kershaw    (\cref{fig:kershaws}) & 47K & 7 & 16M & 6 & 2.6M\\
  146 pebble  (\cref{fig:ns_cases}a) & 62K & 7 & 21M & 6 & 3.5M\\
  1568 pebble (\cref{fig:ns_cases}b) & 524K & 7 & 180M & 72 & 2.5M \\
  67 pebble   (\cref{fig:ns_cases}c) & 122K & 7 & 42M & 18 & 2.3M\\
\bottomrule
\end{tabular}
\end{table}

For the pressure-Poisson tests, three flow cases are considered, as depicted
in \cref{fig:ns_cases}. 
The first three cases corresponds to turbulent flow through a cylindrical packed-bed
with 146, 1568, and 67 spherical pebbles.
The 146 and 1568 pebble cases are from Lan and coworkers \cite{lan2021all}.
The 67 pebble case is constructed using an alternate Voronoi cell approach,
and includes chamfers \cite{RegerATH2022}.
As such, the 67 pebble case is a more complex geometry.
The first two bed flows are at Reynolds number $Re_D=5000$, based on sphere
diameter, $D$, while the 67 pebble case is at Reynolds number $Re_D=1460$. Time advancement is based on a two-stage 2nd-order
characteristics timestepper with CFL=4 
($\dt = 2\times10^{-3}$
$\dt = 5\times 10^{-4}$, and
$\dt = 5\times 10^{-5}$ for the 146, 1568, and 67 pebble cases). 
An absolute pressure solver tolerance of $10^{-4}$ is used. A restart at
$t=10$, $t=20$, and $t=10$ convective time units is used for the 146, 1568, and 67 pebble
cases, respectively, to provide an initially turbulent flow.  

In all cases, solver results are collected over 2,000 timesteps. At each step,
the solution is projected onto a space of up to 10 prior solution vectors to
generate a high-quality initial guess, ${\bar \uu}$. Projection is standard
practice in nekRS as it can reduce the initial residual by orders of magnitude
at the cost of just two matrix-vector products in $A$ per step
\cite{fischer1998projection}.
The solver used is GMRES(15) preconditioned with a single pMG V-cycle.

\section{High Order Results}
\label{sec:results}
Here we consider the solver performance results for the test cases
of \cref{sec:cases}. We assign a single MPI rank to each GPU
and denote the number of ranks as $P$. All runs are on Summit.
Each node on Summit consists of 2 IBM POWER9 processors and 6 NVIDIA V100 GPUs.
Each case is preconditioned using pMG with a schedule of $(7,5,3,1)$
for Jacobi-based and $(7,3,1)$ for Schwarz-based Chebyshev smoothing.
While $p=1$ is treated as the coarse grid level for pMG,
the problem sizes still scale linearly with the number of spectral elements.
At the $p=1$ level, boomerAMG is used on the CPU is used with the settings
described in \cref{sec:semfem}.

\subsection{Loss of Symmetry and GMRES versus CG} \label{sec:loss-of-symmetry}
One concern regarding the one-sided V-cycle approach is the loss
of symmetry in the preconditioner, requiring the use of GMRES instead of CG as a KSP solver.
There are, however, several strategies to mitigate this issue.
Using restarted GMRES($m$) bounds the orthogonalization cost
for $n$ degrees of freedom to at most $O(m^2 n)$ per iteration.
Further, an effective preconditioning strategy can reduce this cost further
by ensuring that $k$ is small.
An additional concern using GMRES($m$) is the $O(m)$ all-reduce operations per iteration.
However, by utilizing classical Gram-Schmidt orthogonalization, GMRES($m$) requires only two all-reduce operations per iteration.
Concerns over the loss of orthogonality are avoided by keeping $m$ small (e.g., $m=15$ or $m=30$).
While not used in this current work, Thomas and coworkers demonstrated
that %
post-modern GMRES reduces the number of synchronizations to a single all-reduce
while preserving backward stability \cite{thomas2022post}.

Kershaw (\cref{sec:kershaw}) results utilizing $1^{st}$ Cheb, Jacobi(3,3), (7,5,3,1), pMG as a \emph{symmetric} preconditioner
with CG and GMRES(30) as the KSP solvers are shown in \cref{tab:ksp-comparison}.
Since the Kershaw case described in \cref{sec:kershaw} uses a $10^{-8}$ residual reduction
for the convergence criterion of the solver, it follows that GMRES should yield a lower
iteration count than CG.
This is because CG minimizes the $A$-norm of the error vector, while GMRES minimizes the $L_2$-norm of the error vector \cite{saad2003iterative}.
A more meaningful convergence criteria would be on the $L_2$-norm of the error vector,
however, this would not be known.
As such, the results presented in \cref{tab:ksp-comparison} are meant to illustrate
that the time to apply a single GMRES iteration is similar to CG.
Consider that, for this case, GMRES(30) is employed, thereby \emph{over-estimating}
the cost per iteration for the NS cases wherein GMRES(15) is used.
Lastly, the most effective preconditioning strategies utilize Schwarz-based Chebyshev smoothing.
The ASM and RAS methods considered herein are \emph{asymmetric}, and thus are not suitable for use with CG.
\begin{table}
\centering
\caption{\small \label{tab:ksp-comparison}Comparison of CG and GMRES(30) for the Kershaw cases using
$1^{st}$-Cheb, Jacobi(3,3), (7,5,3,1) as preconditioner.
}
\begin{tabular}{rrrrrrr}
\toprule
 $\varepsilon$ & CG iter. & Time per CG iter. &  GMRES iter. & Time per GMRES iter. \\
\midrule
  1            &   20 &  $1.26\times 10^{-2}$ &   9 &  $1.27\times 10^{-2}$ \\
  0.3          &  286 &  $1.38\times 10^{-2}$ & 123 &  $1.50\times 10^{-2}$ \\
  0.05         & 1000 &  $1.39\times 10^{-2}$ & 474 &  $1.52\times 10^{-2}$ \\
\bottomrule
\end{tabular}
\end{table}

\subsection{Kershaw Results}\label{sec:kershaw-results}
\begin{figure}
  \def\svgwidth{\textwidth}
  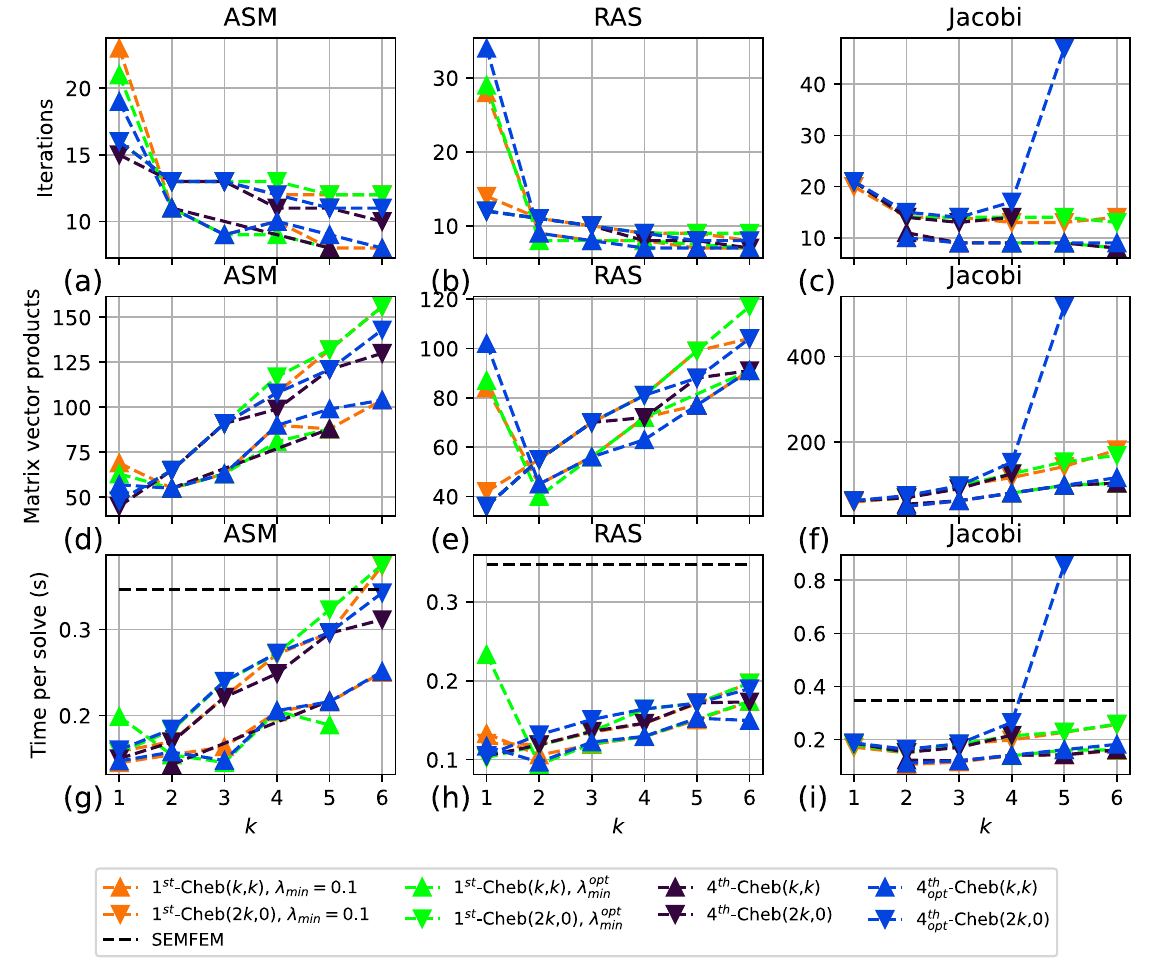
  \caption{\small \label{fig:kershaw-eps-1} Kershaw results, $\varepsilon=1$.}
\end{figure}

\begin{figure}
  \centering
  \def\svgwidth{\textwidth}
  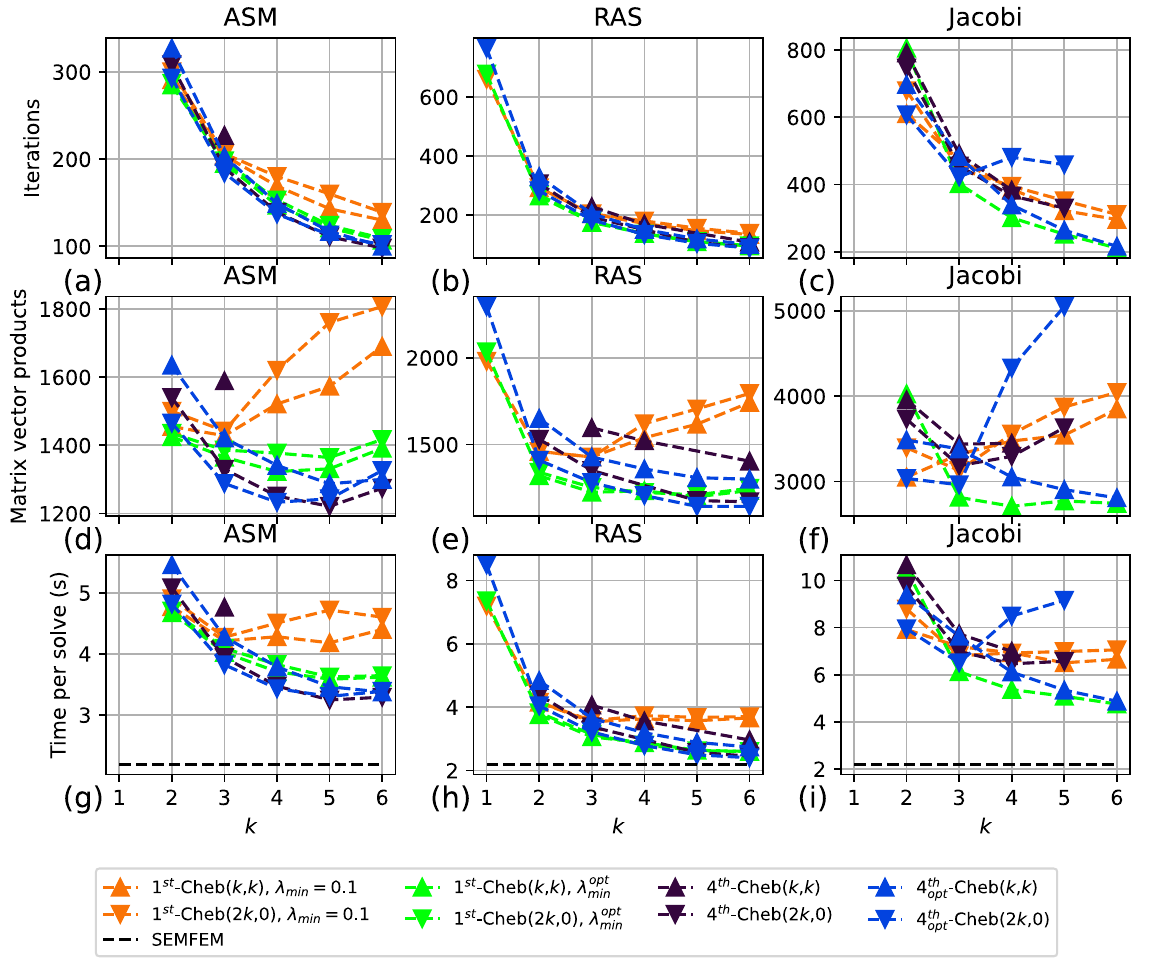
  \caption{\small \label{fig:kershaw-eps-0.05} Kershaw results, $\varepsilon=0.05$.}
\end{figure}

Results for the Kershaw case for $\varepsilon=1,0.05$ are shown in
\cref{fig:kershaw-eps-1} and \cref{fig:kershaw-eps-0.05}
\footnote{
  Results for Kershaw with $\varepsilon=0.3$ are shown in \cref{fig:kershaw-eps-0.3}.
}
.
To mitigate the effects of system noise, the reported solve times are based on minimum time to solution over 50 trials.
The reported number of matrix-vector products are for the \emph{finest grid}, $p=7$.

When $\varepsilon=1$, the time to solution is minimized by utilizing
a symmetric V-cycle with relatively low-order Chebyshev-accelerated RAS smoothing.
For this case, SEMFEM is a comparatively poor method.
While the iteration count is decreased with respect to increasing order,
the overall cost of applying the heavier smoother translates to a higher cost per solve.
This can especially be observed in the increased work requirement
in terms of matrix-vector products, \cref{fig:kershaw-eps-1}d-f.
At larger scales, however, the scalability of the AMG coarse grid solve may dominate
the cost of the preconditioner, and thus anything to reduce the number of coarse grid solves
may prove beneficial.
In this scenario one should expect the multigrid approximation property constant \cref{eq:C-def}
to be quite low for the case with no geometric deformation.
In this regime, the theoretical prediction in \cref{fig:critical-c}, states that the convergence is improved
using a symmetric $(k,k)$ V-cycle with as opposed to the one-sided $(2k,0)$ V-cycle.
Benchmarking demonstrates that a single boomerAMG V-cycle iteration for $p=1$ on the CPU
is nearly 12 times the cost of a matrix-vector product.
For larger cases, wherein the number of AMG levels needed for a single
boomerAMG V-cycle increases, the relative cost of the coarse grid solve
will increase relative to the cost of a matrix-vector product.

For $\varepsilon = 0.05$, the fastest time to solution
is reached using SEMFEM preconditioning.
The pMG methods are significantly more expensive.
For the 4th and opt. 4th-kind Chebyshev-accelerated RAS schemes considered,
the time to solution is lowered by increasing the order.
Remarkably, the 4th-kind and optimized 4th-kind Chebyshev schemes
are generally equivalent to, if not better than, optimizing the $\lambda_{min}$ parameter for the standard $1^{st}$ Chebyshev scheme.
This allows for increased performance with high order Chebyshev smoothing, without
the requirement of tuning an additional parameter.
\footnote{
  In large-scale fluid mechanics applications, this overhead is easily
  amortized over the $10^4$--$10^6$ timesteps required.
  Reasonably good values of $\lambda_{min}^{opt}$ do not depend on the RHS,
  allowing this to be part of the setup cost of the preconditioner.
}
Further, with the extreme geometric deformation, the multigrid approximation property constant \cref{eq:C-def}
is expected to be quite large for this case.
The results in \cref{fig:critical-c}
predict that the one-sided approach yields a better convergence rate.
The results confirm this theoretical expectation.

A summary of the results for the Kershaw case is shown in \cref{tab:fastest_solver_nekrs}.
This table reports the solver yielding the lowest time to solution (in seconds), $T_S$,
the iteration count, and the speedup over the time to solution of the default nekRS solver, $T_D$.
The ratio between the time spent
doing coarse grid solves for the default solver, $\left(T_{crs}\right)_D$,
and the time spent doing coarse grid solves for the fastest solver, $\left(T_{crs}\right)_S$,
is reported.
The default nekRS solver is $1^{st}$-Cheb, ASM(3,3),(7,3,1).
For the Kershaw case, these tables demonstrate that optimizing $\lambda_{min}$
in the 1st-kind Chebyshev scheme greatly improves the solver performance.
While $\varepsilon=1,0.3$ do not benefit from the use of a one-sided V-cycle, the use
of the one-sided V-cycle, in conjunction with the optimized 4th-kind Chebyshev
smoother, is able to increase the solver speedup relative to the default solver
by another 13\% for $\varepsilon=0.05$. A 75\% speedup is achieved over the default solver
for $\varepsilon=1,0.05$. For $\varepsilon=0.3$, a much more modest 35\% is achieved.

\subsection{Navier-Stokes}

\begin{figure}
  \def\svgwidth{\textwidth}
  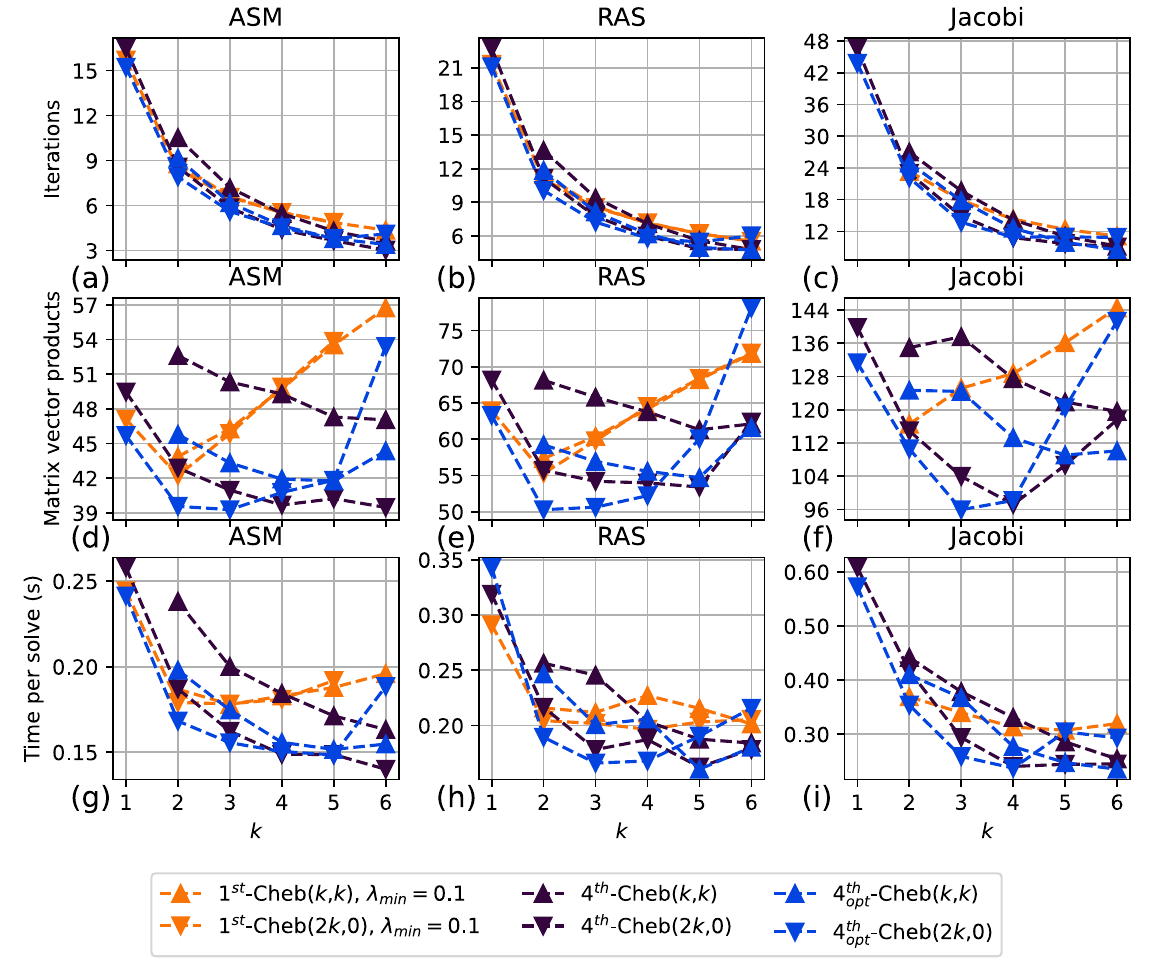
  \caption{\small \label{fig:pb1568-results} 1568 pebble results.}
\end{figure}

\begin{figure}
  \def\svgwidth{\textwidth}
  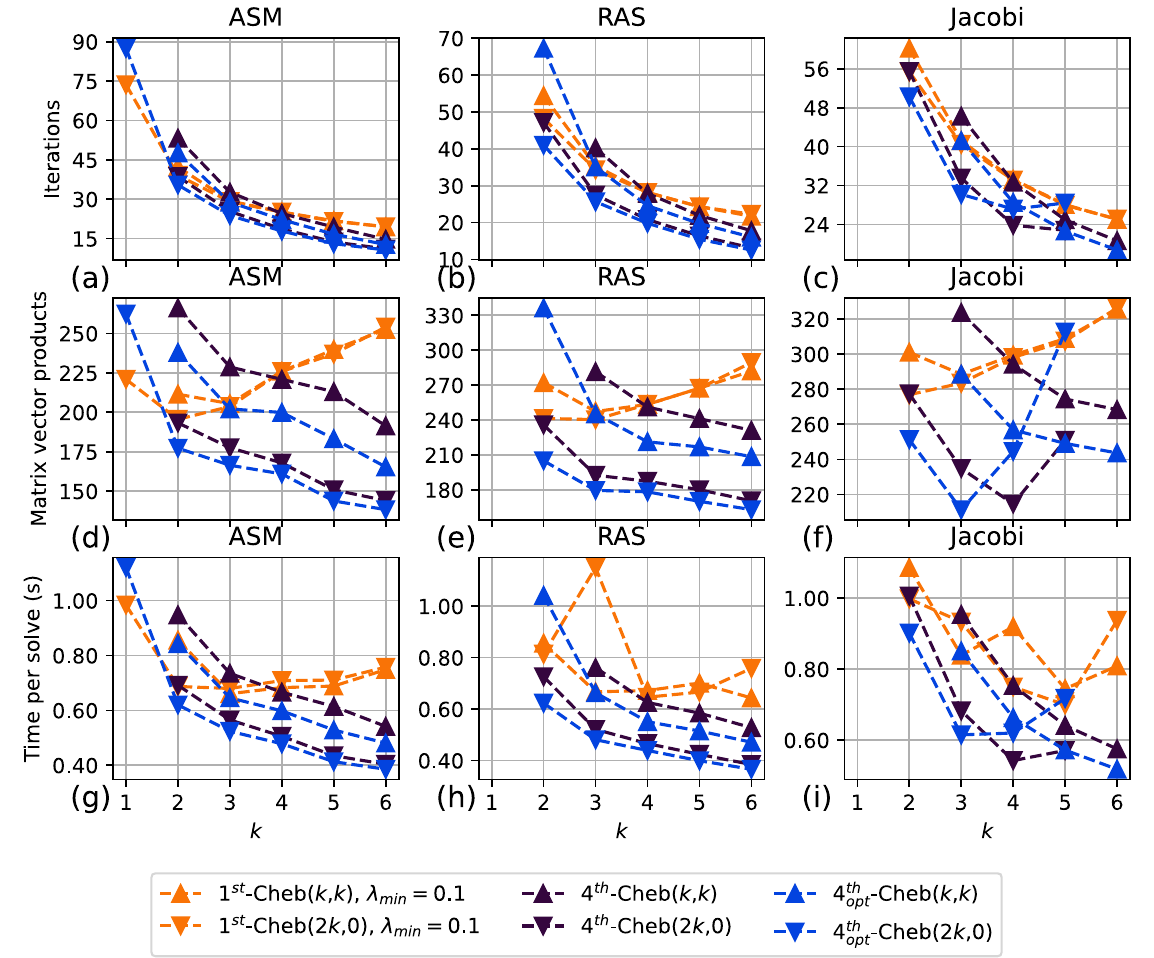
  \caption{\small \label{fig:pb67-results} 67 pebble results.}
\end{figure}

Results for the 1568 and 67 pebble cases are shown
in \cref{fig:pb1568-results} and \cref{fig:pb67-results},
respectively
\footnote{
  Results for the 146 pebble case are shown in \cref{fig:pb146-results}.
}
.
Since the 4th and optimized 4th-kind Chebyshev smoothers are comparable
to the 1st-kind Chebyshev smoother with optimized $\lambda_{min}$ as shown in \cref{sec:kershaw-results},
this smoother is omitted from these cases.

The lowest time to solution for the 1568 pebble case is achieved using
$4^{th}$-Cheb, ASM(12,0), \cref{tab:fastest_solver_nekrs}.
\cref{fig:pb1568-results}d,g shows that the number of matrix-vector products remains nearly
constant with respect to the order, yielding a lower
time to solution by minimizing the coarse grid cost.
The performance of the 1st-kind Chebyshev smoother, however, plateaus
around $k=3, \tilde k = 6$, especially for the Schwarz-based smoothers, see \cref{fig:pb1568-results}g,h.
Although some improvement is observed through the use of the one-sided V-cycle,
an additional benefit is observed by using the alternate Chebyshev smoothers
from Lottes's work \cite{lottes2022optimal}.
Without the added benefit of the one-sided V-cycle, the fastest solver
for this case yields a 17\% speedup over the default. %
However, enabling this one-sided V-cycle approach further increases the solver performance
to a 27\% speedup over the default solver (\cref{tab:fastest_solver_nekrs}).

The 67 pebble case is distinct from the other two pebble meshes.
The 146 and 1568 pebble cases are meshed using an all-hex meshing
strategy developed by Lan and coworkers \cite{lan2021all}.
The 67 pebble case includes chamfers, and is meshed using an alternate Voronoi cell approach \cite{RegerATH2022}.
The resulting mesh proves to be much more challenging for the pMG
based strategies.
Nevertheless, pMG with one-sided Chebyshev smoothing
prove promising for this case.
The time to solution is minimized at relatively high orders
with RAS-based Chebyshev.
This is further improved through the use of the one-sided V-cycle.
The \emph{total work} required per solve is reduced at higher orders
for the Schwarz-based approaches.
Through tuning the pMG parameters, a significant speedup of 81\% is achieved
over the default solver.
In this case, the fastest pMG preconditioner is $4^{th}_{opt}$-Cheb, RAS(12,0).
\begin{table}\small
\centering
\caption{\label{tab:fastest_solver_nekrs}Solver configuration with the fastest time to solution.}
\begin{tabular}{llrrrr}
\toprule
             Case            &                                     Fastest Solver                 & $T_S$          &  Iterations &  $\dfrac{T_D}{T_S}$ & $\dfrac{\left(T_{crs}\right)_D}{\left(T_{crs}\right)_S}$\\
\midrule
Kershaw ($\varepsilon=1   $) & $1^{st}$-Cheb, $\lambda_{min}^{opt}$, RAS(2,2) &           0.09 &        8 &     1.75                & 1.13     \\
Kershaw ($\varepsilon=0.3 $) & $1^{st}$-Cheb, $\lambda_{min}^{opt}$, RAS(5,5) &           0.67 &       28 &     1.35                & 1.79     \\
Kershaw ($\varepsilon=0.05$) &        $4^{th}_{opt}$-Cheb, RAS(12,0)          &           2.40 &       88 &     1.75         & 2.31     \\
\hline
146 pebble                        &         $4^{th}_{opt}$-Cheb, RAS(4,4)     &           0.15 &        5.3 &     1.17         & 1.21  \\
67 pebble                         &        $4^{th}_{opt}$-Cheb, RAS(12,0)     &           0.37 &       12.5 &     1.81  & 2.41  \\
1568 pebble                       &            $4^{th}$-Cheb, ASM(12,0)       &           0.14 &        3 &     1.27    & 2.13  \\
\bottomrule
\end{tabular}
\end{table}

\section{Conclusions}
\label{sec:conclusions}
In this work, we explore the use of one-sided V-cycle using Chebyshev smoothing
as a preconditioner, both in a 2D finite difference example, as well as 
a pMG preconditioner for the pressure Poisson equation arising from the spectral
element discretization of the NS equation.
Further, we explore the efficacy of novel Chebyshev smoothers based on the work
of Lottes \cite{lottes2022optimal} and demonstrate their improvement
over the 1st-kind Chebyshev smoothers.
This improvement enables increased solver performance, especially at high
Chebyshev degrees.
The benefit to this approach is further decreasing the
coarse grid cost.

The authors plan on applying the Lottes's novel Chebyshev smoothers
to AMG solvers, such as BoomerAMG \cite{yang2002boomeramg} and Trilinos/MueLu \cite{prokopenko2014muelu}.
The application of this work are two fold.
First, improved AMG solvers will benefit the low-order SEMFEM
preconditioning strategy for high-order finite elements.
Secondly, improvements in AMG solvers will benefit the overall
solver community.
The authors also plan on applying the ideas developed herein to larger problems.

\section*{Acknowledgments}
This research is supported by the Exascale Computing Project (17-SC-20-SC),
a collaborative effort of two U.S. Department of Energy organizations (Office of Science and the National Nuclear
Security Administration) responsible for the planning and preparation of a capable exascale ecosystem,
including software, applications, hardware,
advanced system engineering and early testbed
platforms, in support of the nation's exascale computing imperative.
This research also used resources of the Oak Ridge Leadership Computing Facility at Oak Ridge National
Laboratory, which is supported by the Office of
Science of the U.S. Department of Energy under
Contract DE-AC05-00OR22725.

The authors thank James Lottes for his insightful comments and suggestions,
in addition to his reading of the manuscript.
The authors thank YuHsiang Lan,
David Alan Reger, and Haomin Yuan
for providing visualizations and mesh files.
\bibliographystyle{siamplain}
\bibliography{references}

\section{Supplementary Material}
\subsection{Coefficients for the Optimized 4th-kind Chebyshev Smoother}

Values for $\beta$ for the optimized 4th-kind Chebyshev smoother developed by Lottes are tabulate in \cref{tab:beta_coef}.

\begin{table}[h]
\centering
\caption{\small \label{tab:beta_coef}Tabulated values of $\beta$ for the optimized 4th-kind Chebyshev smoother.}
\begin{tabular}{ll}
  \toprule
  $k$ & $\beta_i^{(k)}$ \\
  \midrule
  1 & 1.12500000000000 \\
  \hline
  2 & 1.02387287570313 \\
    & 1.26408905371085 \\
  \hline
  3 & 1.00842544782028 \\
    & 1.08867839208730 \\
    & 1.33753125909618 \\
  \hline
  4 & 1.00391310427285 \\
    & 1.04035811188593 \\
    & 1.14863498546254 \\
    & 1.38268869241000 \\
  \hline
  5 & 1.00212930146164 \\
    & 1.02173711549260 \\
    & 1.07872433192603 \\
    & 1.19810065292663 \\
    & 1.41322542791682 \\
  \hline
  6 & 1.00128517255940 \\
    & 1.01304293035233 \\
    & 1.04678215124113 \\
    & 1.11616489419675 \\
    & 1.23829020218444 \\
    & 1.43524297106744 \\
  \hline
  7 & 1.00083464397912 \\
    & 1.00843949430122 \\
    & 1.03008707768713 \\
    & 1.07408384092003 \\
    & 1.15036186707366 \\
    & 1.27116474046139 \\
    & 1.45186658649364 \\
  \hline
  8 & 1.00057246631197 \\
    & 1.00577427662415 \\
    & 1.02050187922941 \\
    & 1.05019803444565 \\
    & 1.10115572984941 \\
    & 1.18086042806856 \\
    & 1.29838585382576 \\
    & 1.46486073151099 \\
  \hline
  9 & 1.00040960072832 \\
    & 1.00412439506106 \\
    & 1.01460212148266 \\
    & 1.03561113626671 \\
    & 1.07139972529194 \\
    & 1.12688273710962 \\
    & 1.20785219140729 \\
    & 1.32121930716746 \\
    & 1.47529642820699 \vspace{0.35cm} \\
  \bottomrule
\end{tabular}
\begin{tabular}{ll}
  \toprule
  $k$ & $\beta_i^{(k)}$ \\
  \midrule
  10 & 1.00030312229652 \\
     & 1.00304840660796 \\
     & 1.01077022715387 \\
     & 1.02619011597640 \\
     & 1.05231724933755 \\
     & 1.09255743207549 \\
     & 1.15083376663972 \\
     & 1.23172250870894 \\
     & 1.34060802024460 \\
     & 1.48386124407011 \\
  \hline
  11 & 1.00023058595209 \\
     & 1.00231675024028 \\
     & 1.00817245396304 \\
     & 1.01982986566342 \\
     & 1.03950210235324 \\
     & 1.06965042700541 \\
     & 1.11305754295742 \\
     & 1.17290876275564 \\
     & 1.25288300576792 \\
     & 1.35725579919519 \\
     & 1.49101672564139 \\
  \hline
  12 & 1.00017947200828 \\
     & 1.00180189139619 \\
     & 1.00634861907307 \\
     & 1.01537864566306 \\
     & 1.03056942830760 \\
     & 1.05376019693943 \\
     & 1.08699862592072 \\
     & 1.13259183097913 \\
     & 1.19316273358172 \\
     & 1.27171293675110 \\
     & 1.37169337969799 \\
     & 1.49708418575562 \\
  \hline
  13 & 1.00014241921559 \\
     & 1.00142906932629 \\
     & 1.00503028986298 \\
     & 1.01216910518495 \\
     & 1.02414874342792 \\
     & 1.04238158880820 \\
     & 1.06842008128700 \\
     & 1.10399010936759 \\
     & 1.15102748242645 \\
     & 1.21171811910125 \\
     & 1.28854264865128 \\
     & 1.38432619380991 \\
     & 1.50229418757368 \\
  \bottomrule
\end{tabular}
\begin{tabular}{ll}
  \toprule
  $k$ & $\beta_i^{(k)}$ \\
  \midrule
  14 & 1.00011490538261 \\
     & 1.00115246376914 \\
     & 1.00405357333264 \\
     & 1.00979590573153 \\
     & 1.01941300472994 \\
     & 1.03401425035436 \\
     & 1.05480599606629 \\
     & 1.08311420301813 \\
     & 1.12040891660892 \\
     & 1.16833095655446 \\
     & 1.22872122288238 \\
     & 1.30365305707817 \\
     & 1.39546814053678 \\
     & 1.50681646209583 \\
  \hline
  15 & 1.00009404750752 \\
     & 1.00094291696343 \\
     & 1.00331449056444 \\
     & 1.00800294833816 \\
     & 1.01584236259140 \\
     & 1.02772083317705 \\
     & 1.04459535422831 \\
     & 1.06750761206125 \\
     & 1.09760092545889 \\
     & 1.13613855366157 \\
     & 1.18452361426236 \\
     & 1.24432087304475 \\
     & 1.31728069083392 \\
     & 1.40536543893560 \\
     & 1.51077872501845 \\
  \hline
  16 & 1.00007794828179 \\
     & 1.00078126847253 \\
     & 1.00274487974401 \\
     & 1.00662291017015 \\
     & 1.01309858836971 \\
     & 1.02289448329337 \\
     & 1.03678321409983 \\
     & 1.05559875719896 \\
     & 1.08024848405560 \\
     & 1.11172607131497 \\
     & 1.15112543431072 \\
     & 1.19965584614973 \\
     & 1.25865841744946 \\
     & 1.32962412656664 \\
     & 1.41421360695576 \\
     & 1.51427891730346 \\
     & \\ %
  \bottomrule
\end{tabular}
\end{table}

\subsection{Geometric $p$ Multigrid for Pressure Poisson}

Mesh quality metrics for the high-order cases are shown in \cref{table:mesh-metrics}.
\begin{table}
\centering
\caption{
  Mesh quality metrics for cases from \cref{fig:kershaws} and \cref{fig:ns_cases}.
  \label{table:mesh-metrics}
}
\begin{tabular}{lcc}
  \toprule
  Case Name 
    & Scaled Jacobian (min/max/avg) 
    & Aspect Ratio (min/max/avg)\\
  \midrule%
  K. ($\varepsilon=1$)
    & 1 / 1 / 1
    & 1 / 1 / 1 \\
  K. ($\varepsilon=0.3$)
    & 0.316 / 1 / 0.841
    & 1.08 / 20.1 / 4.64 \\
  K. ($\varepsilon=0.05$)
    & $1.86\times 10^{-2}$ / 1 / 0.733
    & 1.1 / 162 / 21.7 \\
  \hline
  146 pebble
    & $4.31\times 10^{-2}$ / .977 / .419
    & 1.07 / 56.9 / 7.14 \\
  1568 pebble
    & $2.59\times 10^{-2}$ / .99 / .371
    & 1.12 / 108 / 12.6 \\
  67 pebble
    & $5.97\times 10^{-3}$ / .970 / .38
    & 1.17 / 204 / 13.2 \\
\bottomrule
\end{tabular}
\end{table}

\subsubsection{Poisson}

\begin{figure}[h]
  \def\svgwidth{\columnwidth}
  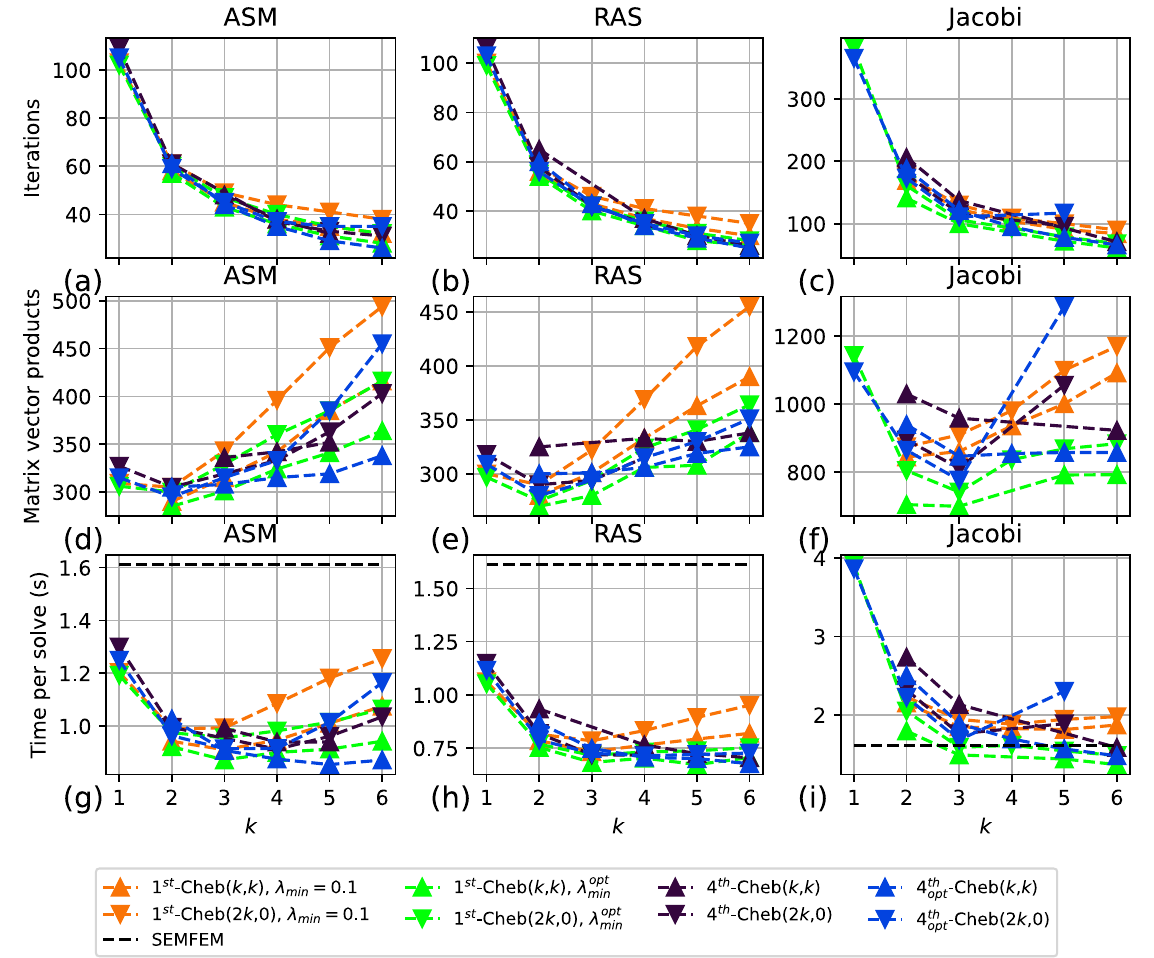
  \caption{\small \label{fig:kershaw-eps-0.3} Kershaw results, $\varepsilon=0.3$.}
\end{figure}

Kershaw results for $\varepsilon=0.3$ are shown in \cref{fig:kershaw-eps-0.3}.
When $\varepsilon=0.3$, however, higher-order Chebyshev smoothing can yield a
lower time to solution. As shown in \cref{fig:kershaw-eps-0.3}d, the time to solution
for the $4^{th}$ and $4^{th}_{opt}$ Chebyshev schemes tend to improve with high orders,
even up to $k=6$ for the symmetric V-cycle (or $\tilde k=12$ for the one-sided V-cycle).
This demonstrates a major improvement over the standard $1^{st}$ Chebyshev scheme,
which tends to yield a minimum time to solution at $k=3$.
Similar to the $\varepsilon=1$ case, SEMFEM is not the preconditioner yielding the
fastest time to solution. However, this approach becomes comparable to Jacobi-based Chebyshev smoothing.

\begin{figure}[h]
  \def\svgwidth{\columnwidth}
  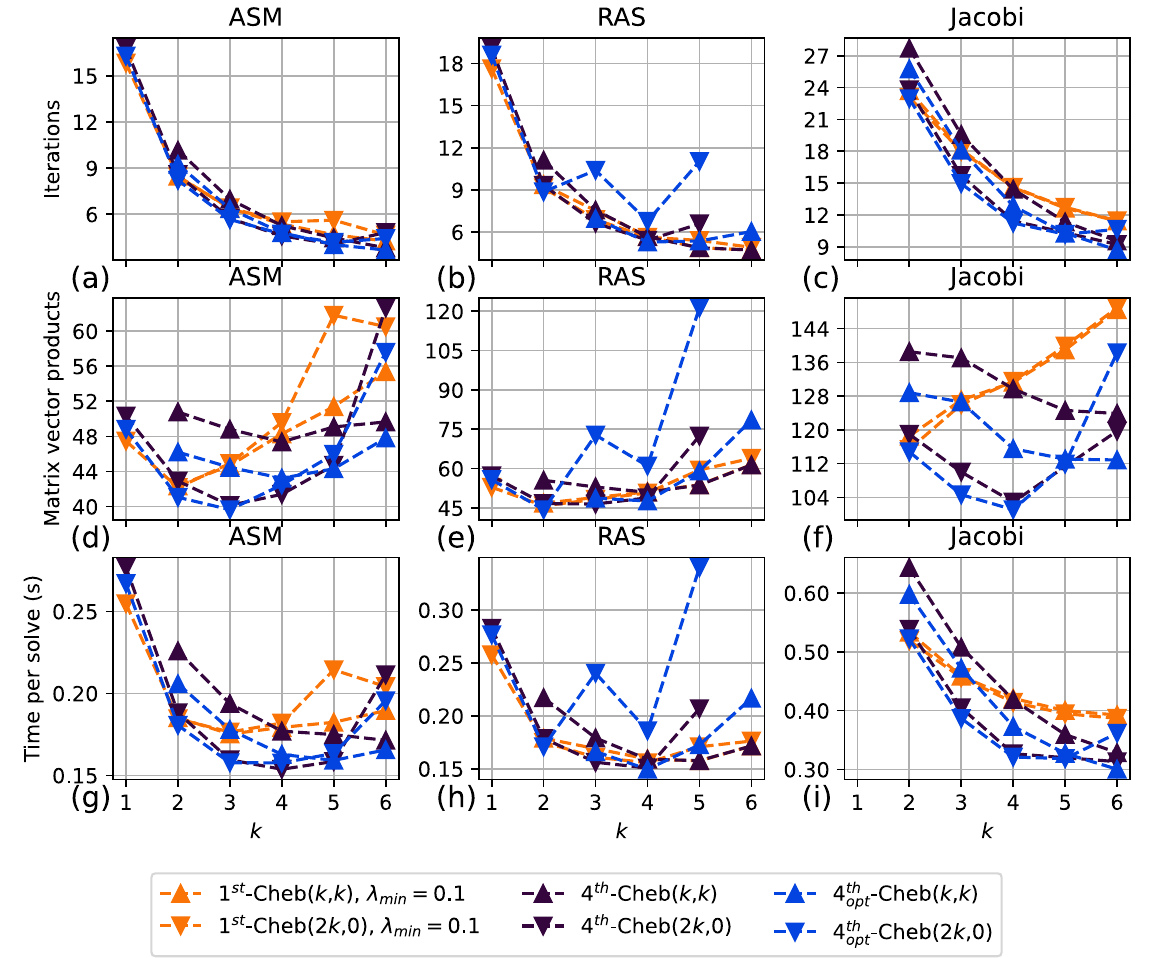
  \caption{\small \label{fig:pb146-results} 146 pebble results.}
\end{figure}

\subsubsection{Navier-Stokes}

pMG results for the 146 pebble Navier-Stokes case mentioned
in \cref{sec:navier-stokes} are shown in \cref{fig:pb146-results}.

\end{document}